\documentclass[aos]{imsart} 

\RequirePackage{amsthm,amsfonts,amssymb}

\RequirePackage[colorlinks,citecolor=blue,urlcolor=blue]{hyperref}
\RequirePackage{graphicx}

\usepackage[reqno]{amsmath}

\usepackage{dsfont}
\usepackage[usenames,dvipsnames]{color}
\usepackage{mathrsfs}
\usepackage{placeins}
\usepackage{graphicx}
\usepackage{tikz}
\usetikzlibrary{arrows}
\usepackage{url}
\usepackage{comment}
\usepackage{multicol}
\usepackage{stmaryrd}
\usepackage{float}
\restylefloat{table}

\usepackage{caption}
\usepackage{subcaption}
\usepackage{natbib}
\usepackage{enumitem}
\setlist[enumerate,1]{label=(\alph*)}

\usepackage{aligned-overset}
\usepackage{bbm}
\usepackage{longtable} 
\usepackage{booktabs}
\usepackage{cleveref}
\crefname{assumptionletter}{Assumption}{Assumption}
\Crefname{assumptionletter}{Assumption}{Assumption}
\usepackage{multirow}

\usepackage{bibunits}

\defaultbibliographystyle{imsart-nameyear} 
\defaultbibliography{Informationskriterien_Bib}

\usepackage{etoolbox} 
\apptocmd{\sloppy}{\hbadness 10000\relax}{}{}

\crefname{equation}{}{} 
\Crefname{equation}{}{} 
\theoremstyle{plain} 
\newtheorem{theorem}{Theorem}[section]
\newtheorem{lemma}[theorem]{Lemma}
\newtheorem{proposition}[theorem]{Proposition}
\newtheorem{corollary}[theorem]{Corollary}

\newtheorem{assumptionTinner}{Assumption}

\Crefname{assumptionTinner}{Assumption}{Assumptions}
\crefname{assumptionTinner}{Assumption}{Assumptions}
\newenvironment{assumptionTa}[1]
  {\begin{assumptionTinner}\label{#1}}
  {\end{assumptionTinner}}

\newtheorem{assumptionMinner}{Model}

\Crefname{assumptionMinner}{Model}{Models}
\crefname{assumptionMinner}{Model}{Models}
\newenvironment{assumptionMa}[1]
  {\begin{assumptionMinner}\label{#1}}
  {\end{assumptionMinner}}

\newtheorem{assumptionDinner}{Assumption}

\Crefname{assumptionDinner}{Assumption}{Assumptions}
\crefname{assumptionDinner}{Assumption}{Assumptions}

\Crefname{assumptionD}{Assumption D}{Assumption D}
\crefname{assumptionD}{Assumption D}{Assumption D}

\Crefname{assumptionT}{Assumption T}{Assumption T}
\crefname{assumptionT}{Assumption T}{Assumption T}

\newlength{\widestop}
\newcommand{\centeredop}[1]{\mathrel{\makebox[\widestop]{$#1$}}}

\newcommand{\thistheoremname}{}
\newtheorem*{genericthm*}{\thistheoremname}
\newenvironment{namedthm*}[1]
  {\renewcommand{\thistheoremname}{#1}%
   \begin{genericthm*}}
  {\end{genericthm*}}

\theoremstyle{definition}
\newtheorem{definition}[theorem]{Definition}

\newtheorem{example}[theorem]{Example}

\newtheorem{remark}[theorem]{Remark}

\def\tablename{\footnotesize \textsc{Table}}

\newcommand{\bthe}{\begin{theorem}}
\newcommand{\ethe}{\end{theorem}}

\newcommand{\ben}{\begin{enumerate}}
\newcommand{\een}{\end{enumerate}}

\newcommand{\bit}{\begin{itemize}}
\newcommand{\eit}{\end{itemize}}

\newcommand{\beq}{\begin{equation}}
\newcommand{\eeq}{\end{equation}}

\newcommand{\ble}{\begin{lemma}}
\newcommand{\ele}{\end{lemma}}

\newcommand{\bde}{\begin{definition}\rm}
\newcommand{\ede}{\halmos\end{definition}}

\newcommand{\bco}{\begin{corollary}}
\newcommand{\eco}{\end{corollary}}

\newcommand{\bpr}{\begin{proposition}}
\newcommand{\epr}{\end{proposition}}

\newcommand{\brem}{\begin{remark}\rm}
\newcommand{\erem}{\end{remark}}

\newcommand{\bproof}{\begin{proof}}
\newcommand{\eproof}{\end{proof}}

\newcommand{\bexam}{\begin{example}\rm}
\newcommand{\eexam}{\end{example}}

\newcommand{\bfi}{\begin{fig}}
\newcommand{\efi}{\end{fig}}

\newcommand{\btab}{\begin{tab}}
\newcommand{\etab}{\end{tab}}

\newcommand{\beao}{\begin{eqnarray*}}
\newcommand{\eeao}{\end{eqnarray*}\noindent}

\newcommand{\balo}{\begin{align*}}
\newcommand{\ealo}{\end{align*}}

\newcommand{\balm}{\begin{align}}
\newcommand{\ealm}{\end{align}\noindent}

\newcommand{\beam}{\begin{eqnarray}}
\newcommand{\eeam}{\end{eqnarray}\noindent}

\newcommand{\barr}{\begin{array}}
\newcommand{\earr}{\end{array}}

\newcommand{\N}{\mathbb{N}}
\renewcommand\P{\mathbb{P}}

\newcommand{\R}{\mathbb{R}}

\def\bX{\boldsymbol X}
\def\bY{\boldsymbol Y}
\def\bZ{\boldsymbol Z}

\def\hats{\widehat{s}_n}
\newcommand{\Tn}[1]{T^{(n)}_{#1}}

\def\bx{\boldsymbol x}

\def\bw{\boldsymbol w}
\def\bv{\boldsymbol v}

\def\bq{\boldsymbol q}

\def\bH{\boldsymbol H}

\def\bTheta{\boldsymbol \Theta}

\def\bSigma{\boldsymbol \Sigma}

\def\bTn{\boldsymbol{T}^{(n)}}
\def\bXn{\bX^{(n)}}

\DeclareMathOperator*{\argmin}{arg\,min}

\newcommand{\vague}{\stackrel{\lower0.2ex\hbox{$\scriptscriptstyle
                    \it{v} $}}{\rightarrow}}
\newcommand{\weak}{\stackrel{\lower0.2ex\hbox{$\scriptscriptstyle
                    \it{w} $}}{\rightarrow}}
\newcommand{\what}{\stackrel{\lower0.2ex\hbox{$\scriptscriptstyle
                    \it{\hat{w}} $}}{\rightarrow}}
\newcommand{\eqdis}{\stackrel{\lower0.2ex\hbox{$\scriptscriptstyle
                    \mathrm{d}$}}{=}}
\newcommand{\distr}{\stackrel{\lower0.2ex\hbox{$\scriptscriptstyle
                    \it{d} $}}{\rightarrow}}

\newcommand{\Rd}{\mathbb{R}^d_+} 
\newcommand{\Rdn}{\mathbb{R}^{d_n}_+} 
\newcommand{\Sd}{\mathbb{S}^{d-1}_+}

\newcommand{\vb}{\vert  \beta \vert }

\newcommand{\ninf}{n \rightarrow \infty} 
\newcommand{\tinf}{t \rightarrow \infty}

\newcommand{\limn}{\lim_{n \rightarrow \infty}}

\newcommand{\Pconv}{\overset{\mathbb{P}}{\longrightarrow}}

\newcommand{\diag}{\mathrm{diag}}

\newcommand{\sstar}{s^*}
\newcommand{\sn}{{s_n}}
\newcommand{\pn}{{q_n}}
\newcommand{\npi}[1]{\pi_{#1}}
\newlength{\dhatheight}

\renewcommand{\hat}{\widehat}

\DeclareMathOperator{\AIC}{AIC}
\DeclareMathOperator{\BIC}{BIC}
\DeclareMathOperator{\QAIC}{QAIC}

\DeclareMathOperator{\MSEIC}{MSEIC}
\DeclareMathOperator{\BICU}{BICU}
\DeclareMathOperator{\BICL}{BICL}

\allowdisplaybreaks 
\definecolor{darkgreen}{RGB}{0,139,0}

\begin{document}
\begin{bibunit}

\begin{frontmatter}
\title{Statistical inference for extremal directions \vspace*{0.2cm} \\ in high-dimensional spaces}
\runtitle{Statistical inference for extremal directions}

\begin{aug}
  \author{\fnms{Lucas} \snm{Butsch}\ead[label=e1]{lucas.butsch@kit.edu} }
    \and
   \author{\fnms{Vicky} \snm{Fasen-Hartmann}\ead[label=e2]{vicky.fasen@kit.edu}\orcid{0000-0002-5758-1999}}
 \address{Institute of Stochastics, Karlsruhe Institute of Technology \\[2mm] \printead[presep={ }]{e1,e2}}


  \runauthor{L. Butsch and  V. Fasen-Hartmann}
\end{aug}

\begin{abstract}
In multivariate extreme value statistics, the first step in understanding the dependence structure of extremes is identifying the directions in which they occur.    
The novelty of this paper is the analysis of high-dimensional extreme value models in which both the model dimension and the number of bias directions go to infinity as the number of observations tends to infinity; we estimate the number of extremal directions. 
To address the curse of dimensionality, we extend and investigate the information criteria (AIC, BICU, BICL, QAIC and MSEIC) from the fixed-dimensional case \citep{butsch2024fasen,meyer_muscle23}, which employ the concept of sparse regular variation that is closely related to multivariate regular variation, for the estimation of the number of extremal directions. 
For all information criteria, we derive sufficient conditions for consistency. Unlike in the fixed-dimensional case, where only the Bayesian information criteria (BICU and BICL) and the QAIC are consistent, the AIC and MSEIC are also consistent in high dimensions under certain model assumptions. We compare the performance of the different information criteria in a simulation study that includes a detailed analysis of the model assumptions and the necessary and sufficient conditions for consistency.

 \end{abstract}



\end{frontmatter}

\section{Introduction}

Multivariate extreme value statistics analyze the probabilities of joint extreme events in multivariate data. For risk assessment, examining marginal tails alone is insufficient; the dependence structure of extremes must also be quantified, including the probabilities of the directions in which joint extremes appear, that are the probabilities of subvectors where all components are large. On the one hand, the number of possible subvectors is large and grows exponentially. On the other hand, when considering extremes, only the most extreme observations are taken into account, resulting in a relatively small effective sample size. Therefore, the curse of dimensionality is not only apparent in models for extremes, but in high-dimensional settings this poses an even greater challenge for statisticians.

The aim of the present paper is to estimate the number of extremal directions in a high-dimensional model by applying the concepts of sparse regular variation and information criteria. This forms the basis for dimension reduction of multivariate extremes. Approaches for reducing the dimension in multivariate extremes are very versatile. Examples include statistical learning methods \citep{chautru,MR3949044,pmlr-v139-jalalzai21a}, Principal Component Analysis (PCA)  \citep{DS:21,drees2025asymptoticbehaviorprincipalcomponent,butsch2025fasen,Sabourin_et_al_2024,AMDS:2025,CT:19,MR4582715,wan2024characterizing}  and variants of $k$-means \citep{JW:20,Fomichov:Ivanovs,AMDS:24,Bernard_2013}. Recent overviews include \citet{C:S:2025} for an in-depth look at the intersection of machine learning and extreme value theory 
and \citet{Engelke:Ivanovs} for probabilistic and statistical aspects of sparse structures in extremes.

Our line of research employs the concept of sparse regular variation as introduced in  \citet{meyer_sparse} and applied in \citet{meyer_muscle23} and \citet{butsch2024fasen}. It is closely related to multivariate regular variation, a classical concept in multivariate extreme value theory \citep{resnick1987,resnick2007,Falk:Buch}. Suppose $\bX^{(n_0)} $ is a $\R^{d_{n_0}}_+$-valued random vector, and suppose that there exists an index $\alpha > 0$ (tail index) and  a random vector  $\bTheta^{(n_0)}$ on the unit sphere $\mathbb{S}_+^{d_{n_0}-1} \coloneqq \{ \bx \in \R^{d_{n_0}}_+: \ \Vert \bx \Vert = 1 \}$
(in the following we use as norm the $L_1$ norm $\Vert \bx \Vert  \coloneqq \sum_{j=1}^d x_j$ for $\bx = (x_1, \ldots, x_d) \in \Rd$)
such that
\begin{equation} \label{RVMeyer}
\P \left( \frac{\bX^{(n_0)}}{t} > r, \frac{\bX^{(n_0)}}{\Vert\bX^{(n_0)}\Vert} \in A \Big| \, \Vert\bX^{(n_0)}\Vert > t \right) \longrightarrow r^{-\alpha} \P(\bTheta^{(n_0)}\in A) ,\quad t \rightarrow \infty,
\end{equation}
for all $r > 0$ and all Borel sets $A \subset \mathbb{S}_+^{d_{n_0}-1}$ with $\P(\bTheta^{(n_0)}\in \partial A) = 0$. Then, the random vector $\bX^{(n_0)}$ is called \textit{multivariate regularly varying }of index $\alpha$.  The {spectral measure} \linebreak $\P(\bTheta^{(n_0)}\in \cdot)$ contains the information about the dependence structure in the extremes of $\bX^{(n_0)}$ and therefore determining the distribution is a particular goal. In high-dimensional datasets where $d_{n_0}$ is large, this can be challenging and computationally intensive.
Therefore, having some prior information about the support of the spectral measure is of particular interest.
For $\beta \subseteq \{1,\ldots,d_{n_0}\}$ and
\begin{align} \label{C_beta}
    C_\beta \coloneqq \{  \bx \in \mathbb{S}_+^{d_{n_0}-1}  : \bx_\beta > \mathbf{0}, \bx_{\beta^c} = \mathbf{0} \}\subseteq \mathbb{S}_+^{d_{n_0}-1},
\end{align}
a positive probability 
 $\P(\bTheta^{(n_0)}\in C_\beta)>0$ reflects that there is an extremal dependence in the coordinates of $\beta$; they are jointly large and build an extremal cluster. 
 However,  a major challenge in estimating the probability
$\P(\bTheta^{(n_0)}\in C_\beta)$ is due to the fact that $C_\beta$ is not necessarily a continuity set and if $\bX^{(n_0)}$ has a continuous distribution then $\P ( {\bX^{(n_0)}}/{t} > r, {\bX^{(n_0)}}/{\Vert\bX^{(n_0)}\Vert} \in C_\beta | \, \Vert\bX^{(n_0)}\Vert > t )=0$ for $\beta \ne \{1, \ldots, d_{n_0}\}$ .
A conclusion that can be drawn from this is
that for $A=C_\beta$  their exist examples where the asymptotic behavior \eqref{RVMeyer} does not hold, and hence, the empirical estimator is inconsistent.
To circumvent this problem,  \citet{damex} developed the support detection algorithm DAMEX (Detecting Anomalies among Multivariate EXtremes) based on truncated $\varepsilon$-cones to generate continuity sets that approximate the sets $C_{\beta}$. A further alternative is the concept of hidden regular variation on a collection of nonstandard subcones in $[0,\infty]^d\backslash{\{0\}}$ in \citet{tawn}.

The approach of this paper to get information about the extremal directions is based on the concept of \textit{sparse regular variation}
going back to \citet{meyer_sparse}.  The main difference between regular variation and sparse regular variation is that the self-normalization \ $\bX^{(n_0)} / \Vert \bX^{(n_0)} \Vert$ in \eqref{RVMeyer} is replaced by the Euclidean projection $\npi{d_{n_0}} (\bX^{(n_0)}/t)$ of $\bX^{(n)}/t$, which is defined as in  \citet{duchi} as $$
 \npi{d_{n_0}}(\bv)= \argmin_{\bw \in \R_+^{d_{n_0}}: \|\bw\|=1} \lVert \bw - \bv \rVert_2^2\in\mathbb{S}^{d_{n_0}-1}_+ \quad \text{ for }\bv\in\mathbb{R}^{d_{n_0}}_+.$$
 The random vector $\npi{d_{n_0}} (\bX^{(n_0)}/t)$ usually has more zero entries than \ $\bX^{(n_0)} / \Vert \bX^{(n_0)} \Vert$ and therefore, is more sparsely populated, advantageous in high-dimensional models which have a sparse dependence structure in the extremes.
To formalize this,  $\bX^{(n_0)}$ is called \textit{sparse regularly varying}, if a $\mathbb{S}_+^{d_{n_0}-1}$-valued random vector $\bZ^{(n_0)}$   and a non degenerate random variable $R$ exist such that
\begin{equation*}
\P \left(  \frac{ \bX^{(n_0)}}{t} > r , \npi{d_{n_0}}\left(\frac{\bX^{(n_0)}}{ t } \right) \in A \; \middle  \vert  \;   \Vert\bX^{(n_0)}\Vert > t \right) \rightarrow \P ( R > r, \bZ^{(n_0)} \in A), \quad \tinf,
\end{equation*}
for all $r > 0$ and all Borel sets $A \subset \mathbb{S}_+^{d_{n_0}-1} $ with  $\P(  \bZ^{(n_0)} \in \partial A) = 0$. 
Note that $R$ is Pareto$(\alpha)$-distributed for an $\alpha > 0$ and models the radial part, whereas the  $\mathbb{S}_+^{d_{n_0}-1}$-valued random vector $\bZ^{(n_0)}$ corresponds to the angular part. The concepts of multivariate regular variation and sparse regular variation are closely connected. 
Under some mild assumptions, a multivariate random vector is regularly varying of index $\alpha$ if and only if it is sparse regularly varying of index $\alpha$  \citep[Theorem 1]{meyer_sparse}.  Additionally, in this case, the set of maximal directions for the regularly varying random vector $\bX^{(n_0)}$ and the sparse regularly varying random vector $\bX^{(n_0)}$  coincide \citep[Theorem 2]{meyer_sparse}. Note that $\beta\subseteq  \{1,\ldots,d_{n_0}\}$ is a maximal direction under regular variation if
\begin{eqnarray*}
    \P(\bTheta^{(n_0)} \in C_\beta) >0 \quad \text{ and } \quad
    \P(\bTheta^{(n_0)} \in C_{\beta'})=0 \quad \text{ for }
    \beta'\supsetneq \beta.
\end{eqnarray*}
Similarly, in the definition of maximal direction for sparse regular variation, the random vector $\bTheta^{(n_0)}$ is replaced by $\bZ^{(n_0)}$. The same theorem in \citet[Theorem 1]{meyer_sparse}
 says that if $\P(\bTheta^{(n_0)} \in C_\beta) >0$ for $\beta \subseteq \{1,\ldots,d_{n_0}\}$, then $\P(\bZ^{(n_0)} \in C_\beta) >0$.  For a  discrete distribution of $\bTheta^{(n_0)}$, we even have that $\bTheta^{(n_0)}=\bZ^{(n_0)}$ almost surely. 
In conclusion, the $2^{d_{n_0}}-1$ probabilities $p^{(n_0)}(\beta):=\P(\bZ^{(n_0)}\in C_\beta)$ for $\beta\subseteq \{1,\ldots,d_{n_0}\}$ contain the essential information about the  support of $\bTheta^{(n_0)}$ and $\bZ^{(n_0)}$. 

Therefore, it is crucial to know all $\beta \subseteq \{1,\ldots,d_{n_0}\}$ with $p^{(n_0)}(\beta)>0$; such a $\beta$ we call \textit{extremal direction} and we denote by
\begin{eqnarray*}
    s^*:=\vert\mathcal{S}^{(n_0)}\vert \quad \text{ with }\quad \mathcal{S}^{(n_0)}:=\{\beta \subseteq \{1,\ldots,d_{n_0}\}:p^{(n_0)}(\beta)>0\}
\end{eqnarray*}
the \textit{number of extremal directions}. A possible estimator for $p^{(n_0)}(\beta)$ is the empirical estimator defined as follows. Let $\bX^{(n_0)},\bX^{(n_0)}_1,\bX^{(n_0)}_2,\ldots$ be a sequence of independent and identically distributed (i.i.d.) $\R^{d_{n_0}}$-valued random vectors and $(k_n)_{n\in\N}$ be a sequence in $\N$ such that $k_n\to\infty$ and $k_n/n\to 0$ as $n\to\infty$. 
 The number $k_n$ reflects the number of extreme observations used for the estimation, and we denote by
$\Vert \bX_{(1,n)}^{(n_0)}\Vert\geq \cdots\geq \Vert \bX_{(n,n)}^{(n_0)}\Vert $ the order statistic of $\Vert \bX_1^{(n_0)}\Vert, \ldots, \Vert \bX_n^{(n_0)}\Vert $. 
Due to  \citet[Proposition 1]{meyer_muscle23} the empirical estimator
\begin{equation} \label{eq:def_Tn}
\frac{T_n^{(n_0)}(C_\beta, k_n)}{k_n}  \coloneqq \frac{1}{k_n} \sum_{j=1}^{k_n} \mathbbm{1}\left\{  \npi{d_{n_0}}(\bX_j^{(n_0)}/ \Vert \bX_{(k_n+1,n)}^{(n_0)} \Vert) \in C_\beta \right\}\Pconv p^{(n_0)}(\beta),
\end{equation}
is a weakly consistent estimator for $p^{(n_0)}(\beta)$ for any $\beta \subseteq \{1,\ldots,d_{n_0}\}$. However,  under quite general assumptions   \citep[Proposition 2]{meyer_muscle23}, the estimated set of extremal directions
    \begin{equation*}
        \widehat{\mathcal{S}}_n^{(n_0)}\coloneqq \{ \beta\subseteq \{1,\ldots,d_{n_0}\}:  T_n^{(n_0)}(C_{\beta}, k_n)  > 0 \} \quad \text{ with }\quad  \widehat s_n^{(n_0)}:=\vert \widehat{\mathcal{S}}_n^{(n_0)}\vert 
    \end{equation*}
    satisfies 
    \begin{equation} \label{sungleichung}
        \lim_{n \rightarrow \infty} \P(\mathcal{S}^{(n_0)}\subseteq 
            \widehat{\mathcal{S}}_n^{(n_0)})=1,
\end{equation}
so that $\widehat s_n^{(n_0)}$  has the tendency to overestimate $s^*$; it is a biased estimator for  $s^*$. On the one hand, 
for $n$ large, the equality $T_n^{(n_0)}(C_\beta,k_n)=0$
implies that $\beta$ is not an extremal direction. But on the other hand,  if $T_n^{(n_0)}(C_\beta,k_n)>0$ then $\beta$ might be an \textit{extremal direction} or a, so called,  \textit {bias direction}. A key challenge is identifying these different types of directions.  To circumvent the overestimation, \citet{meyer_muscle23} introduced an Akaike information criterion (AIC), the consistency of which was investigated in \citet{butsch2024fasen}. Furthermore, \citet{butsch2024fasen}
developed additional information criteria including the quasi-Akaike information criterion ($\QAIC$), the mean-squared error information criterion ($\MSEIC$) and the Bayesian information criteria ($\BICU$ and $\BICL$), and derived consistency results for these  criteria. In particular, the Bayesian information criteria and the QAIC are consistent information criteria.


The novel aspect of this paper is the estimation of the number of extremal directions in a high-dimensional model. By high-dimensional we mean, that the random vector $\bX^{(n)}$ is an $\R^{d_n}$-valued random vector where $d_n\to\infty$ as $n\to\infty$.  However, for any $n\in\N$, the random vector $\bX^{(n)}$ is still sparsely regularly varying. In this model, additional noise is introduced due to the high dimensions, resulting in an even higher bias in the estimation. As in the context of high-dimensional covariance matrix estimation  \citep{Bai:Choi:Fujikoshi:2018,butsch2025fasen, Bai:Silverstein:2010}, we must assume a certain structure in our model in order to define the parameters properly.
Note that we use the notation $\mathcal{P}_{d}$ to represent the power set of  $\{1,\ldots,d\}$ excluding the empty set. 

\begin{assumptionMa} ~ \label{ModelM}
Let for any  $n\in\N$ the $\R^{d_n}$-valued random vector $\bX^{(n)}$ be sparse regulary varying with index $\alpha$ and $d_n\to\infty$ as $n\to\infty$.
Suppose there exists a $n_0\in\N$ such that for $n\geq n_0$ the following holds:
  \begin{enumerate}
        \item[(M1)] ${\displaystyle \lim_{t\to\infty}\frac{\P(\Vert \bX^{(n)}\Vert >t)}{\P(\Vert \bX^{(n_0)}\Vert >t)}=1}$.
        \item[(M2)] 
   $ p^{(n)}(\beta)=p^{(n_0)}(\beta) $ for all $\beta \in \mathcal{P}_{d_{n_0}}$
     { and } $p^{(n)}(\beta)=0$ for all $\beta\in \mathcal{P}_{d_{n}}\backslash \mathcal{P}_{d_{n_0}}. 
$
    \item[(M3)] $\bX^{(n_0)}$ has $s^*$ extremal directions and for $\beta_j \in \mathcal{P}_{d_{n_0}}$ we have the ordering 
\begin{eqnarray*} 
    p^{(n_0)}(\beta_1)\geq p^{(n_0)}(\beta_2)\geq \ldots\geq 
     p^{(n_0)}(\beta_{s^*})>0=p^{(n_0)}(\beta_{s^*+1})=\ldots=p^{(n_0)}(\beta_{2^{d_{n_0}}-1}).
\end{eqnarray*}
\end{enumerate}
\end{assumptionMa}
One conclusion is that, for any $n\geq n_0$,
the random vector $\bX^{(n)}$ also has the  $s^*$ extremal directions, namely the directions $\beta_1,\ldots,\beta_{s^*}$. All other directions are irrelevant.

This paper focuses on estimating $s^*$ in this high-dimensional \Cref{ModelM} by generalizing the AIC, BIC, MSEIC, and QAIC of \cite{meyer_muscle23} and \cite{butsch2024fasen} and investigating their consistency. Unlike the fixed-dimensional case, we can present sufficient criteria for consistency for all information criteria, which are, except for the MSEIC, as well necessary.
A major challenge is that the empirical estimators  ${T_n^{(n)}(C_\beta, k_n)}/{k_n} $ might be inconsistent estimators for $p^{(n_0)}(\beta)$ as an increase in the dimension introduces additional noise in the estimation process, resulting in a bias. 
This phenomenon is well known when estimating the empirical eigenvalues in high-dimensional covariance models; in such models, the empirical eigenvalues are no longer consistent estimators of the true eigenvalues \citep{Bai:Choi:Fujikoshi:2018,butsch2025fasen, Bai:Silverstein:2010}. This is partly due to the fact that the power set of $\{1,\ldots,d_{n}\}$ is countable and  $\mathbb{S}^{d_{n}-1}_+ = \bigcup_{\beta \in \mathcal{P}_{d_{n}}} C_\beta$
but for $n\to \infty$ the power set $\mathcal{P}(\N)$ of $\N$ is no longer countable.

\subsection*{Example}
To conclude this section, let us consider an example of underestimation. We define for $n\geq n_0$ the sequence of random vectors
    $\bX^{(n)}=(
        \bX^{(n_0)\top},
        X_{d_{n_0}+1},
        \ldots,
        X_{d_n})^\top$,
 where $\bX^{(n_0)}$ is sparse regularly varying of index $\alpha$ with extremal directions $\beta_1,\ldots,\beta_{s^*}$ satisfying (M2) and $X_{n_0+1},X_{n_0+2},\ldots$ is a sequence of i.i.d. random variables with 
\begin{eqnarray*}
    \lim_{t\to\infty}\frac{\P(\vert X_j\vert >t)}{\P(\Vert \bX^{(n_0)}\Vert >t)}=0 \quad \text{ for } \quad  j\geq n_0+1.
\end{eqnarray*}
Then the sequence $(\bX^{(n)})_{n\in\N}$ satisfies   \Cref{ModelM}. 
In the special case where  $(X_j)_{j> n_0}$ is an i.i.d. sequence of positive random variables, where all moments are finite, we obtain the following 
asymptotic behaviors for any $m\geq n_0$:
%
\begin{eqnarray*}
    \lim_{n\to\infty}\frac{\P(\Vert \bX^{(n_0)}\Vert >n)}{\P(\Vert \bX^{(m)}\Vert >n)}=1 
    \quad \text{ and } \quad 
    \lim_{n\to\infty}\frac{\P(\Vert \bX^{(n_0)}\Vert >n)}{\P(\Vert \bX^{(n)}\Vert >n)}=0. 
\end{eqnarray*} 
In this example, we cannot expect the empirical estimators ${T_n^{(n)}(C_{\beta_j}, k_n)}/{k_n}$ for $j=1,\ldots,s^*$, to be consistent estimators of $p^{(n_0)}(\beta_j)$ anymore. This conjecture is confirmed in our simulation study in \Cref{sec:simulation}.

 \subsection*{Structure of the paper} 
 The paper is organized as follows. In \Cref{sec:preliminaries}, we properly define bias direction 
 and introduce our model assumptions.  The main results of the paper are derived in \Cref{sec:information_criteria}, where we propose and analyze the following information criteria: the Bayesian information criteria ($\BICU$ and $\BICL$) in \Cref{sec:BIC}, the Akaike information criterion ($\AIC$) in \Cref{sec:AIC}, the Quasi-Akaike information criterion ($\QAIC$) in \Cref{sec:QAIC}  and finally, the mean-squared error information criterion ($\MSEIC$) in \Cref{sec:MSEIC}. Furthermore, we provide sufficient and necessary conditions for the weak consistency of the information criteria. In \Cref{sec:simulation}, we compare the proposed information criteria in a simulation study, where we consider an asymptotically tail independent model in \Cref{sec:asymp_indep} and an asymptotically tail dependent model in \Cref{sec:asymp_dep}.  In addition, the model assumptions and the necessary and sufficient conditions for consistency are analysed in detail. 
 Finally, in \Cref{sec:conclusion} we come to a conclusion. The main proofs of the paper are provided in the appendix. In the supplementary material, additional simulation studies are presented.

\subsection*{Notation} Throughout the paper, we use the following notation.  
 For a vector $\bx=(x_1,\ldots,x_d)^\top \in \R^d$  we write  $\diag(\bx)\in\R^{d\times d}$ for a diagonal matrix with the components of $\bx$ on the diagonal. Furthermore, $\mathbf{0}_d \coloneqq (0, \ldots, 0)^\top \in \R^d$ is the zero vector. Moreover, $\Vert \bx \Vert \coloneqq \Vert \bx \Vert_1$ is the $L_1$-norm for $\bx\in\R^d$ and the unit sphere $\Sd=\{\bx\in\left[0,\infty\right)^d:\,x_1+\cdots+x_d=1\}$ is defined with respect to the $L_1$-norm. 
  By $| a |$ we denote the absolute value of a real number $a$ and by  $|A|$ the cardinality of a set $A$, but the meaning should be clear from the context. 
In addition, $\mathcal{P}(\N)$ is the power set of $\N$ and
 $\mathcal{P}_d$ is the power set of  $\{1,\ldots,d\}$ excluding, in both cases, the empty set. 
  Finally,  $\Pconv$ is the notation for convergence in probability.

\section{Preliminaries} \label{sec:Extreme:direction} \label{sec:preliminaries} \label{sec:bias direction} 
Firstly, in this section, we extend the fixed-dimensional models of \citet{meyer_muscle23} and \citet{butsch2024fasen} to a high-dimensional model and formulate the precise model assumptions. 
Therefore, suppose that \Cref{ModelM} is given, and that $(k_n)_{n\in\N}$ is a sequence in $\N$ such that $k_n\to\infty$ and $k_n/n\to 0$ as $n\to\infty$ where $k_n$ denotes the number of extreme observations used for the estimation. We use the same notation as in the introduction, except that $C_\beta$, as defined in  \eqref{C_beta}, is now a subset of $\mathbb{S}^{d_{n}-1}_+$
and ${T_n^{(n)}(C_\beta, k_n)}$ is defined as in \eqref{eq:def_Tn} using the i.i.d. sequence $\bX_1^{(n)},\ldots,\bX_n^{(n)}$ with distribution $\bX^{(n)}$ from \Cref{ModelM} so that
\begin{eqnarray*}
    T^{(n)}(C_\beta, k_n):=T_n^{(n)}(C_\beta, k_n)=\frac{1}{k_n} \sum_{j=1}^{k_n} \mathbbm{1}\left\{  \npi{d_{n}}(\bX_j^{(n)}/ \Vert \bX_{(k_n+1,n)}^{(n)} \Vert) \in C_\beta \right\}.
\end{eqnarray*}
If $\beta \in \mathcal{P}(\N) \setminus \mathcal{P}_{d_n}$, then we set $T^{(n)}(C_\beta, k_n): = 0$. 
Note that for any $\beta \in \mathcal{P}(\N)$ the number $\Tn{}(\beta, k_n)$ is bounded by $k_n$ and therefore $\limsup_{\ninf} \Tn{}(\beta, k_n) / k_n \le 1$. 
       For $n\in\N$ the \textit{set of observed directions} is defined as
        \begin{align*}
            \widehat{\mathcal{S}}_n \coloneqq \{ \beta \in \mathcal{P}(\N) : \Tn{} (\beta, k_n) > 0\}
        \end{align*}
        and the \textit{number of observed directions} is defined as
        \begin{align*}
             \hats \coloneqq \vert \widehat{\mathcal{S}}_n  \vert.
        \end{align*}
Under mild assumptions, \eqref{sungleichung} gives that, for the fixed-dimensional case, the set of observed directions overestimates the set of extremal directions. Therefore, we call a $\beta\in \widehat{\mathcal{S}}_n\backslash \{\beta_1,\ldots,\beta_{s^*}\}$ a \textit{bias direction} because it is a observed direction which is not an extremal direction and we expect that this overestimation holds as well in the high-dimensional case.
Moreover, we order the values $\Tn{}(\beta, k_n)$ for $\beta \in \mathcal{P}_{d_n}$ by size such that $\Tn{j}(k_n)$ corresponds to the $j$-th largest value of the sequence $(\Tn{}(\beta, k_n))_{\beta \in \mathcal{P}_{d_n} }$ and define the ordered vector
$
    \bTn (k_n) \coloneqq (\Tn{1}(k_n), \ldots, \Tn{\widehat s_n}(k_n))^\top.
$


\begin{remark} $\mbox{}$
\begin{itemize}
  \item[(a)]  Indeed, if $\hats \Pconv r$ as $\ninf$ with $r < \infty$ then it is possible to project the data onto a finite-dimensional space and receive the same consistency results for the information criteria as in \citet{butsch2024fasen}. Therefore, it is sufficient in the following to investigate the case $\hats \Pconv \infty$. This implies that the number of bias directions is infinite, which is not possible in the fixed-dimensional case. However, we will see that this assumption is confirmed by our simulation study in \Cref{sec:simulation}.
  \item[(b)]  A common approach in high-dimensional statistics is to assume that the model dimension $d_n$ grows to infinity at a rate similar to that of the number of observations $n$  \citep{Bai:Choi:Fujikoshi:2018,butsch2025fasen,Bai:Silverstein:2010}. In the context of multivariate extremes, the dimension is the dimension of $ \bTn (k_n)$, which corresponds to $\widehat s_n$, and the number of observations is the number of extremal observations $k_n$. More precisely, we will assume that  $\hats/k_n$ converges to a constant value bounded by $1$ since $\hats \le \sum_{j=1}^{\hats} \Tn{j}(k_n) \le k_n$. 
\end{itemize}
\end{remark}



For the rest of the paper, we will only evaluate the information criteria up to $q_n$, the so-called \textit{number of candidate models}. 
While evaluating the information criteria across all $s=1,\ldots,\hats $ potentially extremal directions would require excessive computational effort, it is reasonable to assume some prior information on the number of candidate models. Note that possible choices for $q_n$ are constants independent of $n$. 
In the following, we will summarize the assumptions of this paper.


\begin{assumptionTa}{asu:T}
{Suppose \Cref{ModelM} and that 
$\Tn{j}(k_n)=\Tn{}(\beta_j,k_n)$ for $j=1,\ldots,s^*$ and $n\in\N$. 
Additionally, we assume the following: \vspace*{-0.2cm}}
    \begin{enumerate}[label=(T\arabic*)] ~ 
 \item     As $\ninf$ holds $\hats \Pconv \infty$ and $\hats/k_n \Pconv c \in [0,1)$. \label{asu:aic_hd_1}
    \item Suppose that  for  $j=1,\ldots, \sstar$, \label{asu:aic_hd_2}
    \begin{align*}
        \frac{\Tn{j}(k_n)}{k_n} \Pconv  p_j\in (0,1).
    \end{align*}
    \item Suppose that there exists a $\mu \ge 1$ such that  \label{asu:aic_hd_3}
    \begin{align*} 
        \frac{1}{\hats} \sum_{j=\sstar +1}^{\hats} \Tn{j}(k_n) \Pconv \mu.
    \end{align*}
     \item   Suppose that there exists a  $q \ge 1$ such that \label{asu:aic_hd_4}
    \begin{align*} 
    \Tn{\sstar+1}(k_n)  \Pconv q \mu.
    \end{align*}     
    \item Suppose that the number of candidate models $\pn$ satisfies  $$\pn/\sqrt{\hats} = o_\P(1).$$ \label{asu:aic_hd_5} 
    \end{enumerate}\label{asu:aic_high_dim2}
    \vspace*{-0.5cm}
\end{assumptionTa}

\begin{remark}~ \label{rem:asu_T}
We would like to provide some comments to justify these assumptions.
    \begin{enumerate}
    \item The model assumptions are in analogy to the properties of the empirical eigenvalues in high-dimensional covariance models \citep{Bai:Choi:Fujikoshi:2018,butsch2025fasen,Bai:Silverstein:2010} if we replace $\Tn{j}(k_n)/k_n$ by the $j$-th largest empirical eigenvalue. Furthermore, note that if $c=0$, the assumptions \ref{asu:aic_hd_2}-\ref{asu:aic_hd_4} hold in the fixed-dimensional case as has been shown in \cite{meyer_muscle23}  with $p_j=p^{(n_0)}(\beta_j)$ under very general assumptions.
    \item Since $\Tn{1}(k_n)\geq \Tn{2}(k_n)\geq\cdots$, a conclusion of \ref{asu:aic_hd_4} and $k_n\to\infty$ as $n\to\infty$,
    is that for any bias direction $\beta_j$, $j\geq s^*+1$, the empirical estimator
    ${\Tn{j}(k_n)}/{k_n}\Pconv 0$,
    which corresponds to $p^{(m)}(\beta_j)$ if $\beta_j\in\mathcal{P}_{d_{m}}$. Again, this mimics the fixed-dimensional case and the consistency of the empirical estimators in the case of bias directions.
    \item Although any bias direction has an estimated probability mass $\Tn{j}(k_n)/k_n$ converging to zero in probability, we have  that the probability mass lying on all bias directions is accumulating. On the one hand, due to \ref{asu:aic_hd_2},
    \begin{eqnarray*}
        \frac{1}{k_n}\sum_{j=s^*+1}^{\widehat s_n}\Tn{j}(k_n)\Pconv \left(1-\sum_{j=1}^{s^*}p_j\right), 
    \end{eqnarray*}
    and on the other hand, due to \ref{asu:aic_hd_1} and \ref{asu:aic_hd_3}, 
    \begin{eqnarray} \label{cmu}
        \frac{1}{k_n}\sum_{j=s^*+1}^{\widehat s_n}\Tn{j}(k_n)\Pconv c\mu \in\left[0,1\right). 
    \end{eqnarray}
   Only  for $c=0$ the probability mass is negligible with $\sum_{j=1}^{s^*}p_j=1=\sum_{j=1}^{s^*}p^{(n_0)}(\beta_j)$. \item We can interpret $c$ as a measure for the additional noise introduced into the model by increasing its dimension. If $c$ is increasing,  not only does the number of bias directions increase, but the total bias in \eqref{cmu} also increases. This implies that
    the empirical estimator $(\Tn{1}(k_n),\ldots,\Tn{s^*}(k_n))/k_n$ is unable to consistently estimate $(p^{(n_0)}(\beta_1),\ldots,p^{(n_0)}(\beta_1))$ because  $\sum_{j=1}^{s^*}p_j<1$; for some extremal direction $\beta_j$, $j\in\{1,\ldots,s^*\}$, the probability $p^{(n_0)}(\beta_j)$ is underestimated.
   
    \item Although by the definition of a bias direction it might be possible that $\Tn{\sstar+1}(k_n) \Pconv \infty$, we exclude this case in \ref{asu:aic_hd_4}, because then it is even more difficult to distinguish, based on the observations, whether an extremal direction or a bias direction is present. We suspect that the consistency results for the information criteria are analogous to those in the fixed-dimensional case of \cite{butsch2024fasen}.  At least we were able to demonstrate that the AIC is not weakly consistent. For the sake of brevity, however, we will skip the details here.
    \end{enumerate}
\end{remark}

\section{Information criteria} \label{sec:information_criteria}
The aim of this paper is to estimate $\sstar$, the number of extremal directions through information criteria. Therefore, we introduce first the $\BICU$ and $\BICL$ in \Cref{sec:BIC}, then the $\AIC$ in \Cref{sec:AIC}, the $\QAIC$ in \Cref{sec:QAIC}, and finally the $\MSEIC$ in \Cref{sec:MSEIC}. Additionally, we derive necessary and sufficient conditions for consistency for these information criteria. The proofs of this section are provided in Appendix \ref{proof:information_criteria}.

\subsection{Bayesian information criterion} \label{sec:BIC}
 The idea of a Bayesian information criterion ($\BIC$) is to determine the Model $s$ with the highest posterior probability  $\P(  s | \bTn(k_n) )$ given the data $\bTn(k_n)$  for $ s = 1, \ldots, q_n$. Note that we evaluate the information criteria only until $q_n$, the number of candidate models, which is much smaller than $\widehat s_n$. In the fixed-dimensional case, the $\BIC$ was derived in \citet[Section 5.1]{butsch2024fasen} by bounding the posterior probability as in \citet{bicgeneralisation}. By exchanging $r$ by $\hats$ in the definition gives our 
 \textit{Bayesian information criterion concerning the upper bound} ($\BICU$) 
\begin{align*}
    \BICU_{k_n}(s) &\coloneqq  - 2 \log (k_n!) + 2  \sum_{j=1}^{\hats} \log(\Tn{j}(k_n)!) - 2\sum_{j=1}^{s} \Tn{j}(k_n) \log \left( \frac{\Tn{j}(k_n)}{ k_n} \right) \nonumber \\
&\qquad  - 2\log \left( \frac{1}{k_n (\hats -s)} \sum_{j=s+1}^{\hats} \Tn{j}(k_n)  \right)  \sum_{j=s+1}^{\hats} \Tn{j}(k_n)  + 2 s \log \left( k_n  \right) \\
&\qquad + s \log \left(\frac{ \hats}{2 \pi ( \hats-s) }  \right)
\end{align*}
for $s = 1, \ldots, q_n$  and the \textit{Bayesian information criterion concerning the lower bound} ($\BICL$) 
\begin{align*}
    \BICL_{k_n}(s) &\coloneqq  - 2 \log (k_n!) + 2  \sum_{j=1}^{\hats} \log(\Tn{j}(k_n)!) - 2\sum_{j=1}^{s} \Tn{j}(k_n) \log \left( \frac{\Tn{j}(k_n)}{ k_n} \right) \nonumber \\
&\qquad  - 2\log \left( \frac{1}{k_n (\hats -s)} \sum_{j=s+1}^{\hats} \Tn{j}(k_n)  \right)  \sum_{j=s+1}^{\hats} \Tn{j}(k_n)    + s \log \left( k_n  \right) \\
&\qquad + s \log \left(\frac{ k_n}{2 \pi \Tn{1}(k_n)}  \right).
\end{align*}
As in the fixed-dimensional case  \citep[Theorem 5.5]{butsch2024fasen}, the  $\BICU$ and $\BICL$ are weakly consistent in the high-dimensional case as formulated in the next theorem. 

\begin{theorem}\label{th:BIC_Consistency}
    Suppose \Cref{asu:T} holds. Then $\BICU$ and $\BICL$ are weakly consistent in the sense of 
    \begin{enumerate} 
        \item $\displaystyle  \limn  \P\Big(   \argmin_{1 \le s < \pn }  \BICU_{k_n}(s)  = \sstar \Big) = 1,$
        \item $\displaystyle \limn  \P\Big(   \argmin_{1 \le s < \pn }  \BICL_{k_n}(s)  = \sstar \Big) = 1.$
    \end{enumerate}
\end{theorem}
In contrast to the results for spiked covariance models in  \citet{Bai:Choi:Fujikoshi:2018} and \citet{butsch2025fasen}, where the $\BIC$ is only weakly consistent in the high-dimensional setting when a ''gap condition'' is satisfied, in our setting, the  $\BICU$ and $\BICL$ are always weakly consistent in contrast to the $\AIC$ presented below. The empirical eigenvalues of the covariance matrix in these papers correspond here to $T_j^{(n)}/k_n$.

\subsection{Akaike information criterion} \label{sec:AIC}
The Akaike information criterion ($\AIC$) for estimating the number $s^*$  of extremal directions in the fixed-dimensional case was developed in \citet{meyer_muscle23} and is motivated by minimizing the expected Kullback-Leibler divergence between the true distribution of $\bTn(k_n)$  and a multinomial distribution. 
The authors assumed that the number of observed directions is constant, denoted by $r$,  and does not depend on $n$. Therefore, in their definition, we replace $r$ by our estimator $\hats$ and define the \textit{Akaike information criterion} ($\AIC$) as
\begin{align*}
    \AIC_{k_n}(s) &\coloneqq  -  \log (k_n!) +   \sum_{j=1}^{\hats} \log(\Tn{j}(k_n)!) - \sum_{j=1}^{s} \Tn{j}(k_n) \log \left( \frac{\Tn{j}(k_n)}{ k_n} \right) \nonumber \\
&\qquad  - \log \left( \frac{1}{k_n (\hats -s)} \sum_{j=s+1}^{\hats} \Tn{j}(k_n)  \right)  \sum_{j=s+1}^{\hats} \Tn{j}(k_n)  + s ,
\end{align*}
for $s = 1, \ldots, q_n$. 
For fixed $d$, the $\AIC$ is not a weakly consistent information criterion as proved in  \citet[Theorem 3.1]{butsch2024fasen}, which is typical for the $\AIC$ in the large sample size and fixed dimensional case  (cf. \citeauthor{ModelSelection}, \citeyear{ModelSelection}, Section 2.8.2 and \citeauthor{claeskens:16}, \citeyear{claeskens:16}, Section 2.2.1).  However, for high-dimensional spiked covariance models, the $\AIC$ is weakly consistent as derived in \citet{Bai:Choi:Fujikoshi:2018} and \citet{butsch2025fasen}. In our setting, we observe a similar phenomenon and are able to give necessary and sufficient conditions for the consistency in the high-dimensional setup.

\begin{theorem} \label{th:AIC_Consistency}
    Suppose \Cref{asu:T} holds. Then, the AIC is weakly consistent in the sense of     
    \begin{align*}
            \limn  \P\Big(   \argmin_{1 \le s < \pn }  \AIC_{k_n}(s)  = \sstar \Big) = 1,
        \end{align*}
    if and only if     
    \begin{align*}
        g_{\AIC}(q,\mu) \coloneqq  q  \big( 1 -  \log \big(q \big)   \big) -1 + \frac{1}{\mu} > 0.
    \end{align*}
    If $g_{\AIC}(q,\mu)<0$ then
    \begin{align*}
        \limn  \P\Big(   \min_{1 \le s < s^* }  \AIC_{k_n}(s) 
        >\AIC_{k_n}(s^*) \Big) &= 1 \quad \text{ and } \quad\\
        \limn  \P\Big(   \min_{s^*+1 \le s < \pn }  \AIC_{k_n}(s)  > \AIC_{k_n}(s^*) \Big) &= 0.
    \end{align*}
\end{theorem}


\begin{remark}~ \label{rem:aic}
\begin{enumerate}
    \item The consistency of the $\AIC$ relies on the value of $g_{\AIC}(q,\mu)$. Note, if $q < e$, then there exists a $\mu_q > 1$ such that $g_{\AIC}(q,\mu) > 0$ for $\mu \in \left[1, \mu_q\right)$, which results in the consistency of the $\AIC$ for theses pairs $(q,\mu)$. In particular, for $\mu = 1$ and $q \in [1,e)$, the $\AIC$ is weakly consistent whereas for $q >e$, $g_{\AIC}(q,\mu)$ is negative for all $\mu \ge 1$ and hence, the $\AIC$ is inconsistent.
In other words, only for $q$ and $\mu$ sufficiently small, the $\AIC$ is consistent (cf. \Cref{fig:Consistency_HD_for_different_q}). Furthermore, it implies that the number of observations of a bias direction must be sufficiently small to ensure consistency. Hence, there is a big gap between the number of observations of extremal directions and bias directions. This is in agreement with the $\AIC$ consistency results for the number of dominant eigenvalues in a spiked covariance model  \citep{Bai:Choi:Fujikoshi:2018,butsch2025fasen}, where one receives consistency if the gap between the dominant eigenvalues and non-dominant eigenvalues is sufficiently large. 
\item Additionally, it can be shown that if Assumption \ref{asu:aic_hd_4} is not fulfilled, i.e., $\Tn{\sstar+1}(k_n) \Pconv \infty$, then the $\AIC$ is not weakly consistent. 
\end{enumerate}
\end{remark}

\subsection{Quasi-Akaike information criterion} \label{sec:QAIC}
Next, we investigate the quasi-Akaike information criterion ($\QAIC$), which was introduced in \citet[Section 3.1]{butsch2024fasen} in the fixed-dimensional case.  The \textit{quasi-Akaike information criterion} ($\QAIC$) is defined as
\begin{eqnarray*}
    \QAIC_{k_n}(s) &\coloneqq & \hats \log(2 \pi) + \hats \log(k_n) +   \sum_{j=1}^s \log \Big(  \frac{\Tn{j}(k_n)}{k_n} \Big)  + (\hats - s) \log(  \widehat{\rho}^s_n )  + \hats +s
\end{eqnarray*}
for $s = 1, \ldots, q_n$, where $$ \widehat{\rho}^s_n \coloneqq \frac{1}{\hats - s} \sum_{j=s+1}^{\hats}  \frac{\Tn{j}(k_n)}{k_n}.$$ We only replaced $r$ with our estimator $\hats$. Instead of minimizing the Kullback-Leibler distance between the true distribution and the multinomial distribution, as for the $\AIC$, it minimizes the Kullback-Leibler distance between the true distribution and the Gaussian distribution. In the fixed-dimensional setting, the $\QAIC$ is consistent \citep[Theorem 3.8]{butsch2024fasen}, and in the high-dimensional case, we obtain consistency under additional model assumptions.

\begin{theorem}\label{th:QAIC_Consistency}
    Suppose \Cref{asu:T} holds. Then the $\QAIC$ is weakly consistent in the sense of         
            \begin{align*}
                 \limn  \P\Big(   \argmin_{1 \le s < \pn }  \QAIC_{k_n}(s)  = \sstar \Big) = 1,
             \end{align*}
             if and only if
                 \begin{align*}
        g_{\QAIC}(q) \coloneqq \log \left( q \right) -  q  +2  > 0.
    \end{align*}
    If $g_{\QAIC}(q)<0$ then
    \begin{align*}
        \limn  \P\Big(   \min_{1 \le s < s^* }  \QAIC_{k_n}(s) 
        >\QAIC_{k_n}(s^*) \Big) &= 1 \quad \text{ and } \quad\\
        \limn  \P\Big(   \min_{s^*+1 \le s < \pn }  \QAIC_{k_n}(s)  > \QAIC_{k_n}(s^*) \Big) &= 0.
    \end{align*}
\end{theorem}

\begin{remark}\label{rem:QAIC}
The condition $ \log \left( q \right) -  q  +2  > 0$ is similar to the condition for the $\AIC$ in \Cref{th:AIC_Consistency}, where $q$ is not allowed to be too large. $g_{\QAIC}(q)$ is negative if $q > 3,146$ (the value is a numerical approximation as the equation cannot be solved analytically) and consequently the $\QAIC$ is then inconsistent. Hence, again, the number of observations of a bias direction must be sufficiently small. Notably, while the $\QAIC$ is weakly consistent, the interval $(1,\, 3,146)$ is independent of $\mu$. 
\end{remark}

\subsection{Mean squared error information criterion} \label{sec:MSEIC}

The basic idea of the AIC is to minimize the Kullback-Leibler distance of the true distribution and a parametric family of distributions. This minimum is approximated by the expected Kullback-Leibler distance of the true distribution and the estimated distribution, as is done for the $\AIC$. Instead of using the Kullback-Leibler distance, it is also possible to use the normalized mean-squared error (MSE), which results in the $\MSEIC$. The exact derivation is described in \citet[Section 4.1]{butsch2024fasen}. Adapting their approach to our setting by replacing $r$ by $\hats$ gives the  \textit{mean squared error information criterion} ($\MSEIC$) 
\begin{eqnarray*}
\MSEIC_{k_n}(s) &\coloneqq & \frac{k_n}{\sum_{l=s+1}^{\hats} \frac{\Tn{l}(k_n)}{k_n(\hats-s)}} \sum_{j=s+1}^{\hats}  \left( \frac{\Tn{j}(k_n)}{k_n} -  \sum_{i=s+1}^{\hats} \frac{\Tn{i}(k_n)}{k_n(\hats-s)} \right)^{ 2} + 2 s,  
\end{eqnarray*}
for $s = 1, \ldots, q_n$. Similarly to the $\AIC$, the $\MSEIC$ is not weakly consistent in the fixed-dimensional case, but it is in the high-dimensional case.

 \begin{theorem}\label{th:MSE_Consistency}
    Suppose \Cref{asu:T}  holds. The $\MSEIC$ is weakly consistent in the sense of  
             \begin{align*}
                   \limn  \P\Big(   \argmin_{1 \le s < \pn }  \MSEIC_{k_n}(s)  = \sstar \Big) = 1,
     \end{align*}
     if
    \begin{align*}
      g_{\MSEIC}(q,\mu) \coloneqq  2 -  (q - 1)^2 \mu > 0.
    \end{align*}
\end{theorem}

\begin{remark}~ \label{rem:MSEIC}  \label{rem:consistency_functions}
\begin{enumerate}
\item In contrast to the weak consistency of the $\AIC$ (\Cref{th:AIC_Consistency}) and the $\QAIC$ (\Cref{th:QAIC_Consistency}), we proved only a sufficient condition for weak consistency of the $\MSEIC$ but not a necessary condition. The reason is that in the proof of \Cref{th:MSE_Consistency}, in equation \Cref{eq:MSEIC_2}, we are only able to give a lower bound. Otherwise, for the derivation of necessary conditions, we require additional assumptions on the asymptotic behavior of  $\frac{1}{\hats} \sum_{j=\sstar +1}^{\hats} \Tn{j}(k_n)^2$. We believe that under these additional assumptions, our condition on  $g_{\MSEIC}(q,\mu)$ is also necessary. But if the $\MSEIC$
is not consistent, then it overestimates the true value $s^*$.
    \item If $q < 1 + \sqrt{2}$, there exists a $\mu_q> 1$ such that $g_{\MSEIC}(q,\mu) > 0$ for all $\mu \in [1,\mu_q)$, which yields the weak consistency of the $\MSEIC$.  In the special case $\mu = 1$, the $\MSEIC$ is consistent for $q \in [1, 1 + \sqrt{2}) \subset [1,e) \subset [1,\, 3,146)$, where $[1, e)$ is the consistency interval of the AIC and $[1,\, 3,146)$ of the QAIC, respectively  (cf. \Cref{rem:aic} and \Cref{rem:QAIC}).

    \item We see this in \Cref{fig:Consistency_HD_for_different_q}, which shows plots of $g_{\AIC}(q,\mu),g_{\QAIC}(q)$ and $g_{\MSEIC}(q,\mu)$ as functions in $q$ for fixed $\mu$; on the left hand side  for $\mu=1$ and on the right hand side  for $\mu=2$, respectively. If a function is positive, the information criterion is consistent. If $g_{\AIC}(q,\mu)$ or $g_{\QAIC}(q)$, respectively is
negative, then the corresponding information criterion overestimates the number of extremal directions. 
\begin{figure}[ht]
     \centering
         \includegraphics[width=0.8\textwidth]{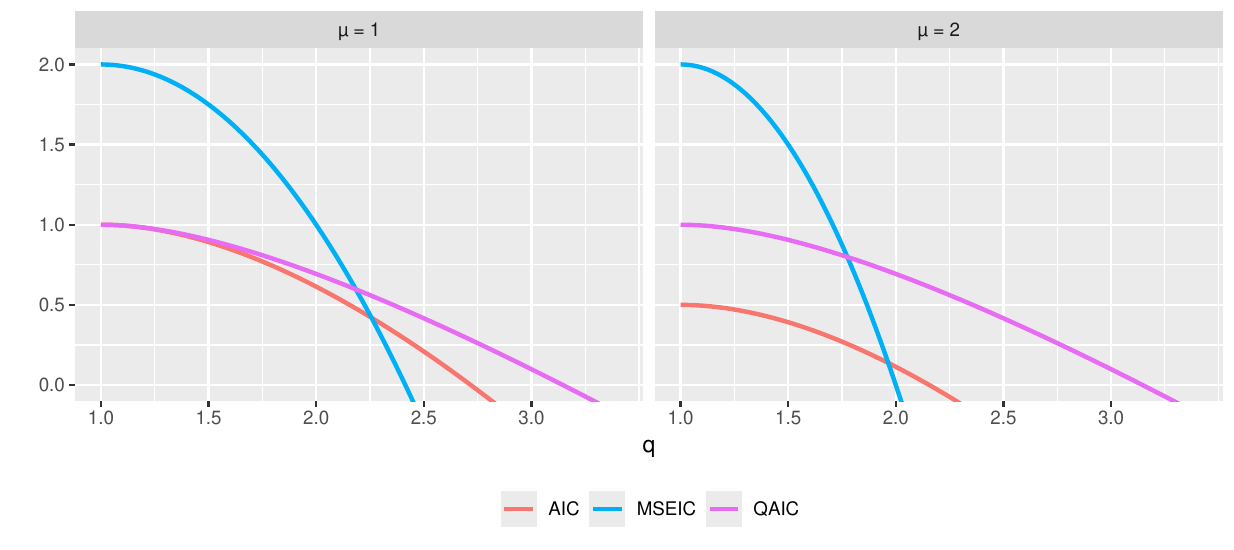}
         \caption{\footnotesize \textit{Plot of   $g_{\AIC}(q,\mu),g_{\QAIC}(q)$ and $g_{\MSEIC}(q,\mu)$, respectively, as functions in $q$ for $\mu =1$ on the left hand side and $\mu = 2$ on the right hand side.}}
    \label{fig:Consistency_HD_for_different_q}
\end{figure} 
\item 
In summary, if $\mu = 1$ and $\mu = 2$, the consistency of the $\MSEIC$ implies the consistency of the $\AIC$ and the $\QAIC$, but opposite conclusions do not hold.  Similarly, the consistency of the $\AIC$ implies the consistency of the $\QAIC$. The consistency intervals for $\mu = 2$ are subsets of the intervals for $\mu =1$.  Note that an increasing $\mu$ reflects that we observe an increase in bias directions. Therefore, it is more difficult for the information criterion to distinguish between a biased direction and an extremal direction and hence, to be weakly consistent.   
 \end{enumerate}
\end{remark}

\section{Simulations} \label{sec:simulation}
In this section, we compare the performance of the proposed information criteria in a simulation study. We consider two different settings, each with $125$ runs. On the one hand, in \Cref{sec:asymp_indep}, we simulate asymptotically independent data, and on the other hand, in \Cref{sec:asymp_dep}, asymptotically dependent data.  Additional simulation studies are provided in the supplementary part. The R-code is available at  \url{https://gitlab.kit.edu/projects/164856} and for the implementation of the projection $\pi$, the approach of \citet{Fast_proj} is used.

\subsection{Asymptotic tail independent model} \label{sec:asymp_indep}

 In the first example, we consider a $d_n$-dimensional random vector 
whose spectral measure has mass on the first $\sstar$ axis. To define the distribution,  we assume that $\bH = (h_{ij})_{1 \leq i,j \leq \sstar} \in \R^{\sstar \times \sstar}$ with  $h_{ij} \overset{ \text{\tiny  i.i.d.}}{\sim} \mathcal{U}((0,1))$ and

\begin{equation*}
\bSigma \coloneqq \diag( h_{11}^{-1/2}, \ldots, h_{dd}^{-1/2}) \cdot \bH^\top \cdot \bH \cdot \diag( h_{11}^{-1/2}, \ldots, h_{\sstar \sstar}^{-1/2}).
\end{equation*}
Note that  $\bSigma_{ii} = 1,\, i = 1, \ldots, \sstar$, and $\bSigma_{ij} < 1, \, i \neq j$. 
Suppose now $\bY=(Y_1,\ldots,Y_{\sstar}) \sim\mathcal{N}_{\sstar} ( \mathbf{0}_{\sstar}, \bSigma)$  under the condition of $\bSigma$ whose components have, by construction, as marginal distribution the standard normal distribution $\Phi$. It is well known that the multivariate normal distribution with correlations smaller than $1$ exhibits pairwise asymptotic independence \citep[Corollary 5.28]{resnick1987}. Now, let 
$\bY_1,\ldots,\bY_{n}, i = 1, \ldots, n$, be an i.i.d. sequence of random vectors with distribution $\bY$ and $N_{i,1}, \ldots, N_{i,d_n-\sstar}$ be the absolute values of an independent sequence of normally distributed random variables. Finally, the i.i.d. random vectors $\bX_i^{(n)} = (X_{i,1}, \ldots, X_{i,d_n})^\top \in \Rdn$, $i = 1, \ldots,n$, are defined as 
\begin{align*}
X_{i,j} \coloneqq
\begin{cases}
  \frac{1}{1 - \Phi ( Y_{i,j})}, & 1 \leq j \leq s^*, \\
  N_{i,j-s^*},                   & s^*+1 \leq j \leq d_n.
\end{cases}
\end{align*}
For our simulation study, we considered $\sstar = 75$ and $d_n  = 200$ leading to $2^{200}-1$ potential directions. 
 The number of candidate models $q_n$ is chosen as $2d_n$, which is small compared to the number of possible dimensions $2^{d_n}-1$.
For the sample sizes $n$ and the number of extremes $k_n$, we used the values $(n,k_n) = (10.000, 550), (25.000, 2.200)$ and $ (50.000,5.000)$.
In these examples, the relative number of extremes $k_n/n$ lies between $5\%$ and $10\%$. We would like to point out a few things about this setup and the relation to our model \Cref{asu:T}: 

\begin{figure}
\centering
 \includegraphics[width=0.9\textwidth]{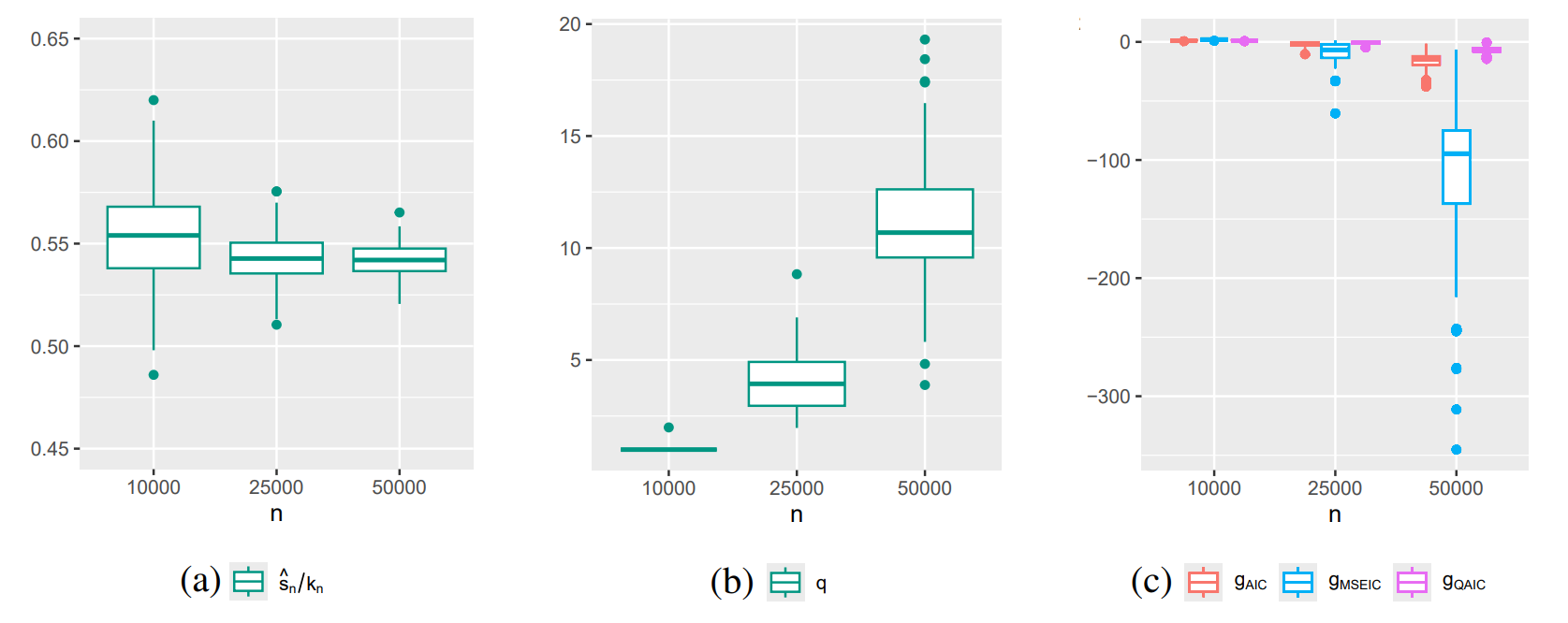}
        
\caption{\footnotesize \textit{
Simulations for the asymptotically independent model with $ s^*  = 75$ and $d_n=200$: Boxplots for the empirical estimator  $\widehat s_n/k_n$ (in (a)), the empirical estimator for $q$ (in (b)), and the empirical estimators for  $ g_{\AIC}(q,\mu),g_{\MSEIC}(q,\mu)$ and $g_{\QAIC}(q)$ (in (c)) are plotted  against the sample sizes $n = 10.000, 25.000$ and $ 50.000$ on the $x$-axis, where
         $ k_{10.000} = 500$, $k_{25.000}=2.200$ and $ k_{50.000}=5.000$, respectively.
}}
\label{fig:asymp_indep_figures}
\end{figure}

\begin{enumerate}
    \item[(i)] Note that we used these values for $k_n$ to achieve nearly constant values of $\hats/k_n$, which can be seen in the boxplot of $\hats/k_n$ in  Figure \ref{fig:asymp_indep_figures}a. On average, $\hats/k_n$ is close to $0,55$ in all simulations and in all simulations $\hats/k_n$ is larger than $0,4$, showing that Assumption \ref{asu:aic_hd_1} is indeed satisfied and that the number of bias directions increases with $n$. Of course, $k_n$ was chosen manually, and the simulations can be improved by a systematic choice of $k_n$
    as is done in the low-dimensional case in \citet{butsch2024fasen}.
    \item[(ii)] We performed simulations as well for $k_n= 0,05\cdot n$, giving similar results and in particular that  $\hats/k_n$ is strictly positive, reflecting that we are in the high-dimensional setup. For the sake of completeness, the results are presented in the supplementary material along with simulations for different dimensions $d_n$. 
    \item[(iii)] Since the average of $\widehat s_n$ is $0,55\cdot k_n$, we have in total on average at most $k_n-(0,55\cdot k_n-s^*)=0,45\cdot k_n+s^*$ observations of any kind of extremal direction and hence, due to the symmetry of our model, every extremal direction has on average at most $(0,45\cdot k_n+s^*)/s^*$
    observations. For $k_n=500,2.200$ and $5.000$, respectively we receive that $(0.45\cdot k_n+s^*)/s^*$ is close to $4$, $14$ and $31$, respectively. Thus, the empirical probabilities of $p_j$ in \ref{asu:aic_hd_2} are at most $(0.45\cdot k_n+s^*)/(s^*k_n)$ and hence, close to $0,008$, $0,0065$ and
    $0,0062$, respectively in the different simulations; they seem to converge to a value strictly less than $1/s^*\approx 0,013$. From this, we see the high-dimensionality of the model again.
    \item[(iv)] In all simulations, the empirical estimator for $\mu$ in \ref{asu:aic_hd_3} was close to 1.
    \item[(v)] From the arguments given in (iii) we know that the empirical mean of $T_{s^*}^{(n)}$ is not higher than $4$, $14$ and $31$ for $n=10.000,25.000$ and $50.000$, respectively, which are also upper bounds for $q$ in \ref{asu:aic_hd_4} if $\mu=1$ (which we observed in the simulations). Figure \ref{fig:asymp_indep_figures}b shows the boxplot for the empirical estimator of $q$ for $n=10.000,25.000$ and $50.000$, respectively. It seems that $q$ is increasing if $n$ increases, which is not surprising, since $k_n$ is increasing as well and thus, $T_{s^*+1}^{(n)}$ is increasing. Note that for $n=10.000$ and $k_n=500$, $q$ is close to 1, since in this setting the observations of both  $T_{s^*+1}^{(n)}$ and $T_{s^*}^{(n)}$ are often 1. An explanation is that $k_n=500$ is quite small in comparison to the number $2^{d_n}-1$ of possible directions.
   \item[(vi)] If $\mu=1$ and $q$ is increasing, the functions $ g_{\AIC}(q,1)$, $g_{\QAIC}(q)$ and $g_{\MSEIC}(q,1)$ are decreasing as we see in \Cref{fig:Consistency_HD_for_different_q}. For $q> 3,146$, we even have $g_{\MSEIC}(q,1)< g_{\AIC}(q,1)<g_{\QAIC}(q)<0$. Due 
   to the empirical estimation of $q$ presented in Figure \ref{fig:asymp_indep_figures}b and the estimation of $\mu$ close to 1, the empirical estimators for $g_{\MSEIC}(q,\mu)$, $g_{\AIC}(q,\mu)$ and $g_{\QAIC}(q)$ plotted in Figure \ref{fig:asymp_indep_figures}c  have this property as well. For $n=50.000$ they are negative, suggesting 
   due to \Cref{th:AIC_Consistency}, \Cref{th:QAIC_Consistency} and \Cref{th:MSE_Consistency}, respectively, that the AIC, QAIC and MSEIC will overestimate the true value $s^*=75$. 
\end{enumerate}

\begin{figure}[ht] 
\centering
\begin{subfigure}{0.49\textwidth}
         \includegraphics[width=\textwidth]{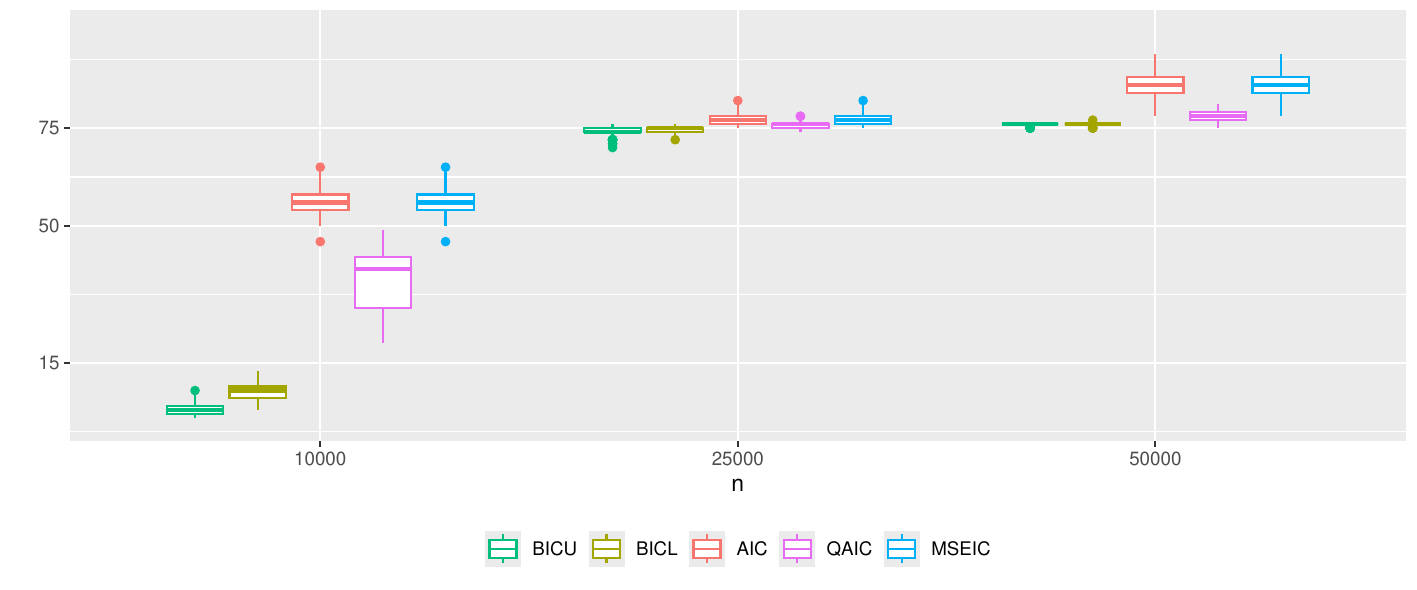}
         \caption{\footnotesize \textit{Number of estimated directions.  
         }}
         \label{fig:asymp_indep_directions}
\end{subfigure}
\hfill
\begin{subfigure}{0.49\textwidth}
         \includegraphics[width=\textwidth]{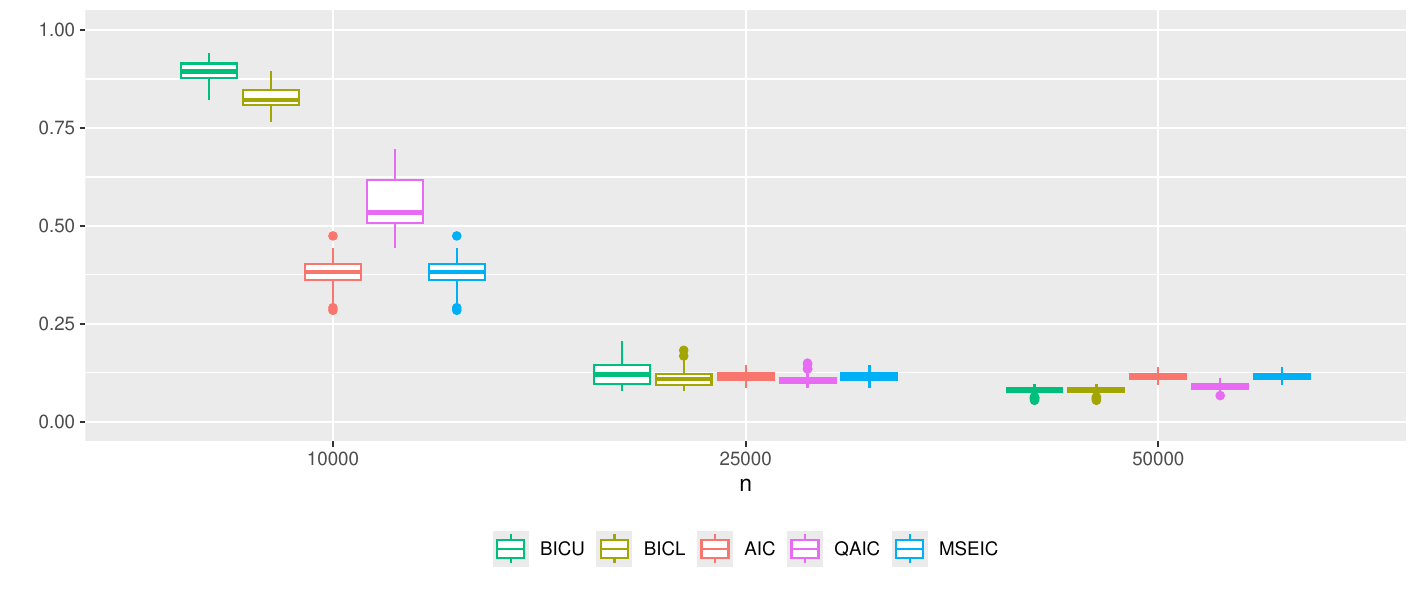}
         \caption{\footnotesize \textit{Hellinger distance.  
         }}
    \label{fig:asymp_indep_hellinger}
\end{subfigure}
         \caption{\footnotesize \textit{
         Simulations for the asymptotically independent model with $ s^*  = 75$ and $d_n=200$: Boxplots for the estimated number of extremal directions (left plot (a)) and the Hellinger distance (right plot (b)), respectively, are plotted against the sample sizes $n = 10.000, 25.000$ and $ 50.000$ on the $x$-axis, where
         $ k_{10.000} = 500$, $k_{25.000}=2.200$ and $ k_{50.000}=5.000$, respectively.}} \label{fig:asymp_indep}
\end{figure}

 Finally, in \Cref{fig:asymp_indep} we see the performance of the different information criteria: On the one hand,  the boxplot of the estimated number of extremal directions $$\widehat s_{\text{IC},n}^*=\argmin_{1\leq s\leq 2d_n}\text{IC}_{k_n}(s)$$ of the different information criteria (IC$_{k_n}$) are presented for $n=10.000,25.000$ and $50.000$ in \Cref{fig:asymp_indep_directions}. On the other hand, in \Cref{fig:asymp_indep_hellinger} we see the boxplots of the Hellinger distance
 \begin{align*}
        H(\widehat{\bq}_{\text{IC},n},\bq)=\frac{1}{\sqrt2} \lVert \widehat{\bq}_{\text{IC},n}- \bq \rVert_2
 \end{align*}
 of the estimated conditional probabilities \begin{equation*} 
    \widehat{\bq}_{\text{IC},n}=  \bigg( \frac{ \Tn{1}(k_n)}{ \sum_{j=1}^{\widehat s_{\text{IC},n}^*} \Tn{j}(k_n)} , \ldots, \frac{ \Tn{\widehat s_{\text{IC},n}^*}(k_n)}{ \sum_{j=1}^{\widehat s_{\text{IC},n}^*} \Tn{j}(k_n)} ,0,\ldots,0\bigg)^{\top}\in[0,1]^{2^{d_n}-1} 
\end{equation*} 
 and the probability vector 
 $\bq:=(\frac{1}{75},\ldots,\frac{1}{75},0,\ldots,0)\in [0,1]^{2^{d_n}-1}$  
 for the different information criteria and sample sizes $n$.
 Note that due to the symmetry of the model, we have $p_{1}=p_2=\ldots=p_{s^*}$ so that $q_j:=p_j/(\sum_{l=1}^{s^*}p_l)=1/s^*=1/75$ for $j=1,\ldots,s^*$
 and $q_j=p_j=0$ for $j>s^*$. This justifies the choice of the true conditional distribution $\bq$ in our asymptotic independent model.
\begin{enumerate}
    \item[(vii)] In both plots of \Cref{fig:asymp_indep}, we see a clear improvement from $n=10.000$ to $n=25.000$ for all information criteria. Where for $n=10.000$ all information criteria are far away from the true value $s^*=75$ and underestimating, for $n=25.000$ the bias and variance in the estimation are quite small for all information criteria. 
    \item[(viii)] If $n=10.000$ we have low observations in any direction (on average not more than 4, cf. (iii)). Therefore, there is not a big difference between   $T_j^{(n)}$ for $j\leq s^*$ and $T_j^{(n)}$ for $j> s^*$, so that for the information criteria it is difficult to detect the gap in $s^*$. If we have a few observations, the Bayesian information criteria perform worse than the other information criteria and select more parsimonious models, due to the heavier penalty. 
   \item[(ix)] In particular, the BICU and BICL improve again from $n=25.000$ to $n=50.000$, which we see in both plots, the estimate of the number of extremal directions in \Cref{fig:asymp_indep_directions} but also in the Hellinger distance in \Cref{fig:asymp_indep_hellinger}: although the variance was already small for $n=25.000$ it is even smaller for $n=50.000$. This shows the consistency of the BICU and BICL derived in \Cref{th:BIC_Consistency}. 
    \item[(x)] In all simulations, the behavior of the AIC and the MSEIC is quite similar. For $n=10.000$ they perform much better than the other information criteria BICU, BICL, and QAIC, although they are not able to estimate the true value $s^*=75$: on average,  AIC and MSEIC estimate $s^*$ as $56,6$, respectively. For $n=25.000$, the average of the AIC and for the MSEIC lies by $76,84$, and for $n=50.000$ it is $86,38$ for both. 
However, an increase in the sample size from $n=25.000$ to $n=50.000$ results
not only in an increase of the bias in the estimation of $s^*$ but also in an increase of the variance, as we see in the growing boxplots in \Cref{fig:asymp_indep_directions}. 
Since the empirical values of $g_{\text{AIC}}(q,\mu)$ and $g_{\text{MSEIC}}(q,\mu)$ in Figure \hyperlink{fig:asymp_indep_figures}{\ref{fig:asymp_indep_figures}c} decrease for $n$ increasing and are strictly smaller than zero for $n=50.000$, the overestimation of the AIC and MSEIC for $n=50.000$ is in accordance with \Cref{th:QAIC_Consistency} and
\Cref{th:MSE_Consistency}, respectively.  The overestimation of the $\AIC$ and the $\MSEIC$ is less pronounced for the Hellinger distance than for the estimation of $s^*$, since the weight assigned to the additional directions is relatively small due to $\mu \approx 1$.  In contrast, for $n=25.000$ the empirical values of $g_{\text{AIC}}(q,\mu)$ and $g_{\text{MSEIC}}(q,\mu)$ are close to zero with a low variance, so that we see only a marginal overestimation in \Cref{fig:asymp_indep_directions}.
\item[(xi)] 
\begin{figure}
\centering
 \includegraphics[width=0.9\textwidth]{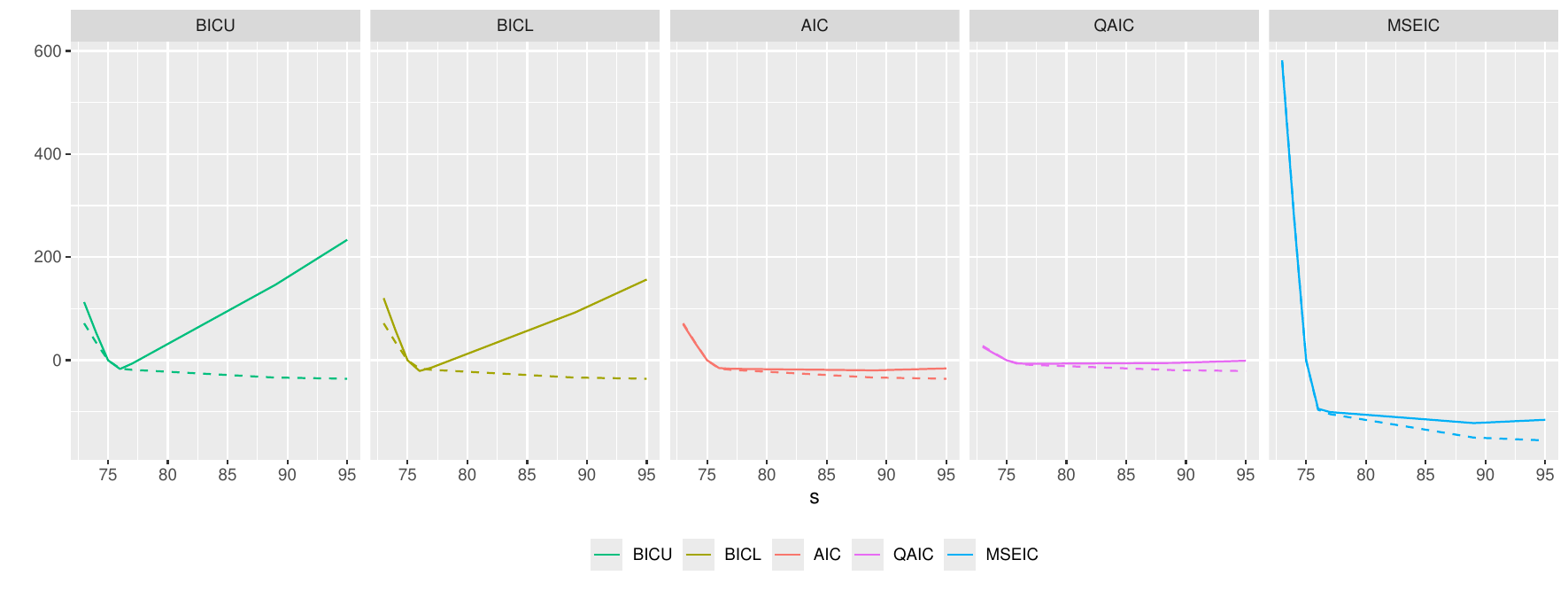}
        
\caption{\footnotesize \textit{
A realization of $\BICU$, $\BICL$, $\AIC$, $\QAIC$ and $\MSEIC$ in the asymptotically independent model with $ s^*  = 75$, $n = 50.000$ and $ d_n=200$: The solid lines indicate $\mbox{IC}_{k_n}(s)-\mbox{IC}_{k_n}(s^*)$ for $73\leq s\leq 95$ and the dashed lines indicate $[(\mbox{IC}_{k_n}(s)-P(s))-(\mbox{IC}_{k_n}(s^*)-P(s^*)]$ (without the penalty term).
}}
\label{fig:asymp_indep_penalty_figures}
\end{figure}
 We realized that in most of the simulations in the high-dimensional setting
\begin{eqnarray*}
    \argmin_{1\leq s\leq \widehat s_n}\mbox{AIC}_{k_n}(s)
     =\argmin_{1\leq s\leq \widehat s_n'}
     \mbox{MSEIC}_{k_n}(s)=\vert \{ j: \Tn{j}(k_n) >1   \}\vert=:\widehat s_n'.        
\end{eqnarray*}
The AIC and MSEIC estimate the number of directions with more than one observation because the penalty is not strong enough for these information criteria. For a better understanding,  we look closer into the sample path behavior of the information criteria: in \Cref{fig:asymp_indep_penalty_figures}, the solid lines indicate  $\mbox{IC}_{k_n}(s)-\mbox{IC}_{k_n}(s^*)$ and the dashed lines
$[(\mbox{IC}_{k_n}(s)-P(s))-(\mbox{IC}_{k_n}(s^*)-P(s^*)]$ for one realization, where $P(s)$ is the penalty term, i.e., for the $\BICU$ the penalty term is $P(s)= 2 s \log \left( k_n  \right)  + s \log \left( \hats/(2 \pi ( \hats-s) )  \right)$, for the $\BICL$ the penalty term is $P(s)= 2 s \log \left( k_n  \right)+ s \log \left( \hats/(2 \pi ( \hats-s))   \right)$, for the $\AIC$ the penalty term is $P(s)=s$, for the $\QAIC$ the penalty term is $P(s)=s$ and finally, for the $\MSEIC$ the penalty term is $P(s)=2s$. In the sample path of  \Cref{fig:asymp_indep_penalty_figures}, the minimum of the AIC and the MSEIC is achieved in 
$\widehat s'_n  = 89. $
Before 76, both information criteria decrease very fast so that the linear penalty $s$ and $2s$, respectively, have no influence. In $s=76$ we see a changepoint in the speed of decrease of both information criteria ($\AIC$, $\QAIC$); the speed of decrease is much lower for $s\geq 76$ than for $s<76$,  but still faster than the slope of the linear penalty term so that the information criteria are still slightly decreasing until $\widehat s'_n  = 89$.
Not surprisingly, this changes after $\widehat s'_n$ again because  for values bigger than $\widehat s_n'$  the likelihood function $\mbox{IC}_{k_n}(s)-P(s)$ of all information criteria is constant (cf. dashed lines in \Cref{fig:asymp_indep_penalty_figures}) and hence,
\begin{eqnarray*}
    \argmin_{1\leq s\leq \widehat s_n}\mbox{IC}_{k_n}(s)
     =\argmin_{1\leq s\leq \widehat s_n'}
     \mbox{IC}_{k_n}(s);        
\end{eqnarray*}
i.e., the minimum is achieved for values less than $\widehat s_n'$. 

      \item[(xii)] We often observe $(1,1,\ldots,1)$  due to the application of the Euclidean projection, so that in total the information criteria estimate 76 instead of 75 extremal directions (cf. the  BIC in  \Cref{fig:asymp_indep_penalty_figures}; the change points in the speed of decrease of the other information criteria lie in 76 instead of 75 in that figure).    
    \item[(xiii)] The QAIC shows overall similar properties to the AIC and MSEIC, but the overestimation and the increase in the variance are not so pronounced (cf. \Cref{fig:asymp_indep}). In \Cref{fig:asymp_indep_penalty_figures} we realize that for $s\leq 76$ the speed of decrease of the $\QAIC$ is much slower than that of $\AIC$ and $\MSEIC$, and hence, the penalty term can compensate for that for $s\geq 76.$
\end{enumerate}

\subsection{Asymptotic dependent model} \label{sec:asymp_dep}
\begin{figure}
\centering
 \includegraphics[width=0.9\textwidth]{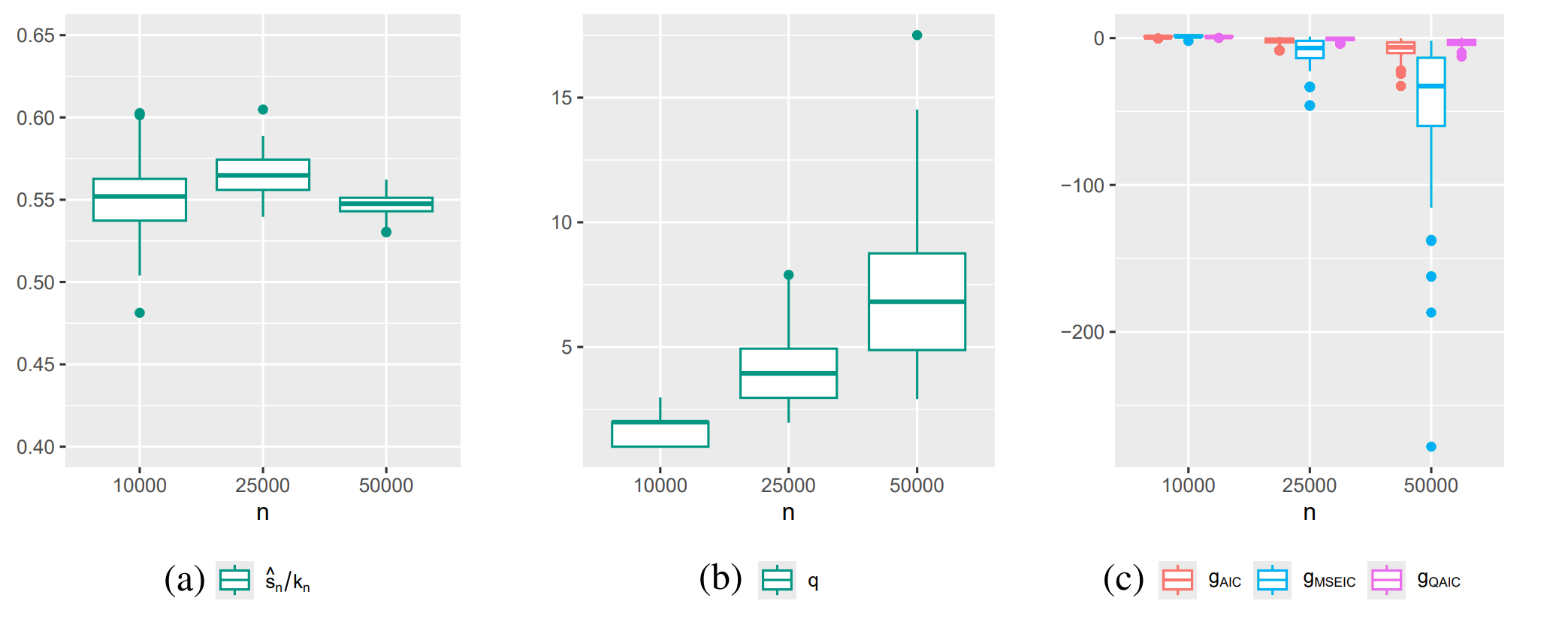}
       
\caption{\footnotesize \textit{Simulations for the asymptotically dependent model with $ s^*  = 50$ and $d_n=300$: The boxplots for  $\widehat s_n/k_n$ (in (a)), the empirical estimator for $q$ (in (b)), and the empirical estimators for  $ g_{\AIC}(q,\mu),g_{\MSEIC}(q,\mu)$ and $g_{\QAIC}(q)$ (in (c)) are plotted against the sample sizes $n = 10.000, 25.000$ and $ 50.000$ on the $x$-axis, where
         $ k_{10.000} = 750$, $k_{25.000}=2.500$ and $ k_{50.000}=5.000$, respectively.}}
\label{fig:asymp_dep_figures}
\end{figure}

\begin{figure}[ht]
\centering
\begin{subfigure}{0.49\textwidth}
         \includegraphics[width=\textwidth]{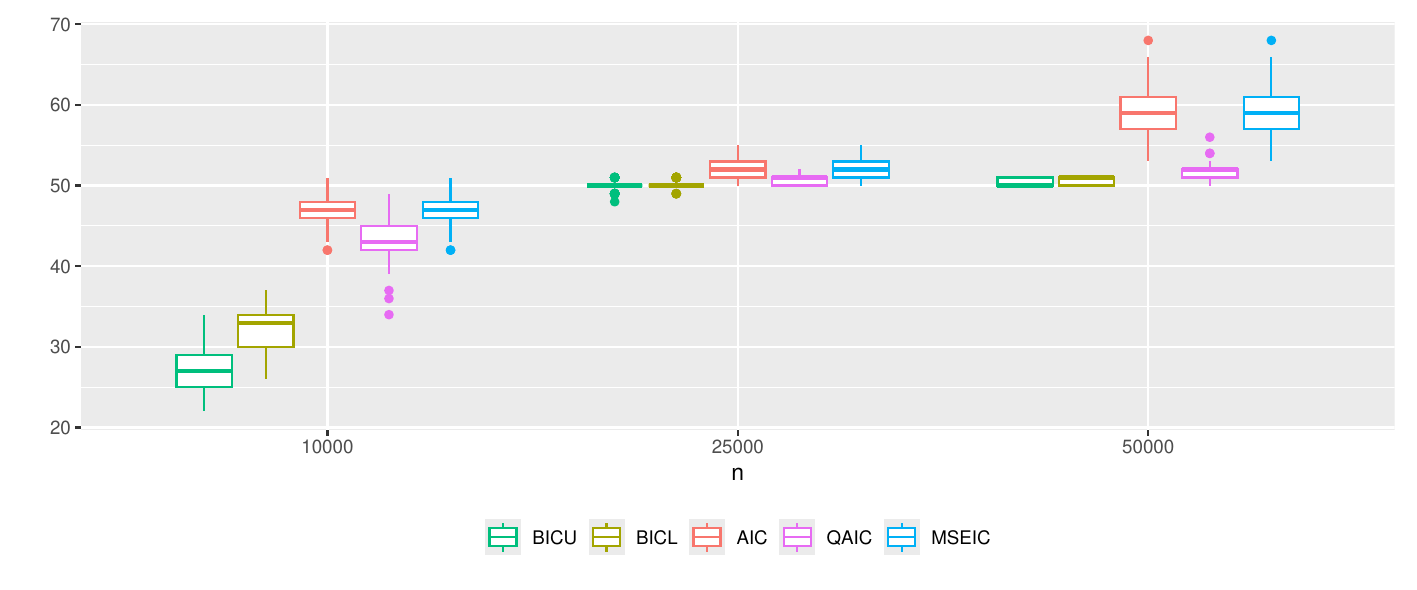}
         \caption{\footnotesize \textit{Number of estimated directions.  
         }}
         \label{fig:asymp_dep_directions}
\end{subfigure}
\hfill
\begin{subfigure}{0.49\textwidth}
         \includegraphics[width=\textwidth]{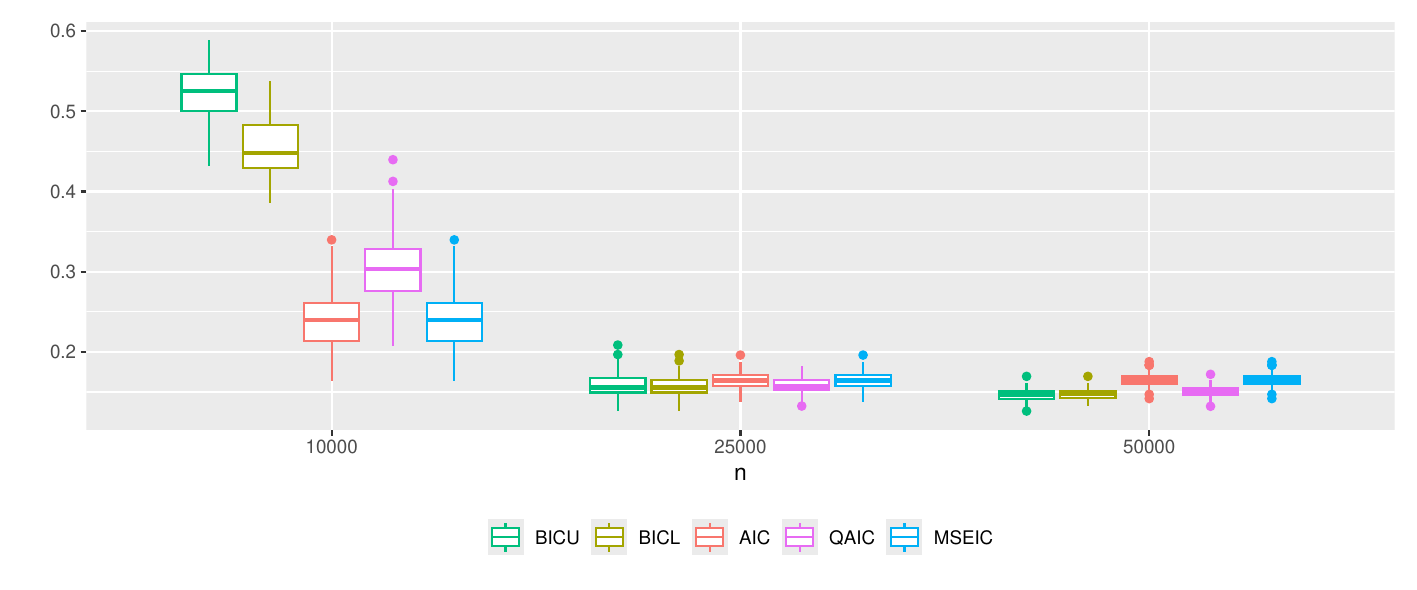}
         \caption{\footnotesize \textit{Hellinger distance.  
         }}
    \label{fig:asymp_dep_hellinger}
\end{subfigure}
         \caption{\footnotesize \textit{Simulations for the asymptotically dependent model with $ s^*  = 50$ and $d_n=300$: Boxplots for estimated number of extremal directions (left plot (a)) and the Hellinger distance (right plot (b)), respectively are plotted  against the sample sizes $n = 10.000, 25.000$ and $ 50.000$ on the $x$-axis, where
         $ k_{10.000} = 750$, $k_{25.000}=2.500$ and $ k_{50.000}=5.000$, respectively.
         }}
         \label{fig:asymp_dep}
\end{figure}

\begin{figure}
\centering
 \includegraphics[width=0.9\textwidth]{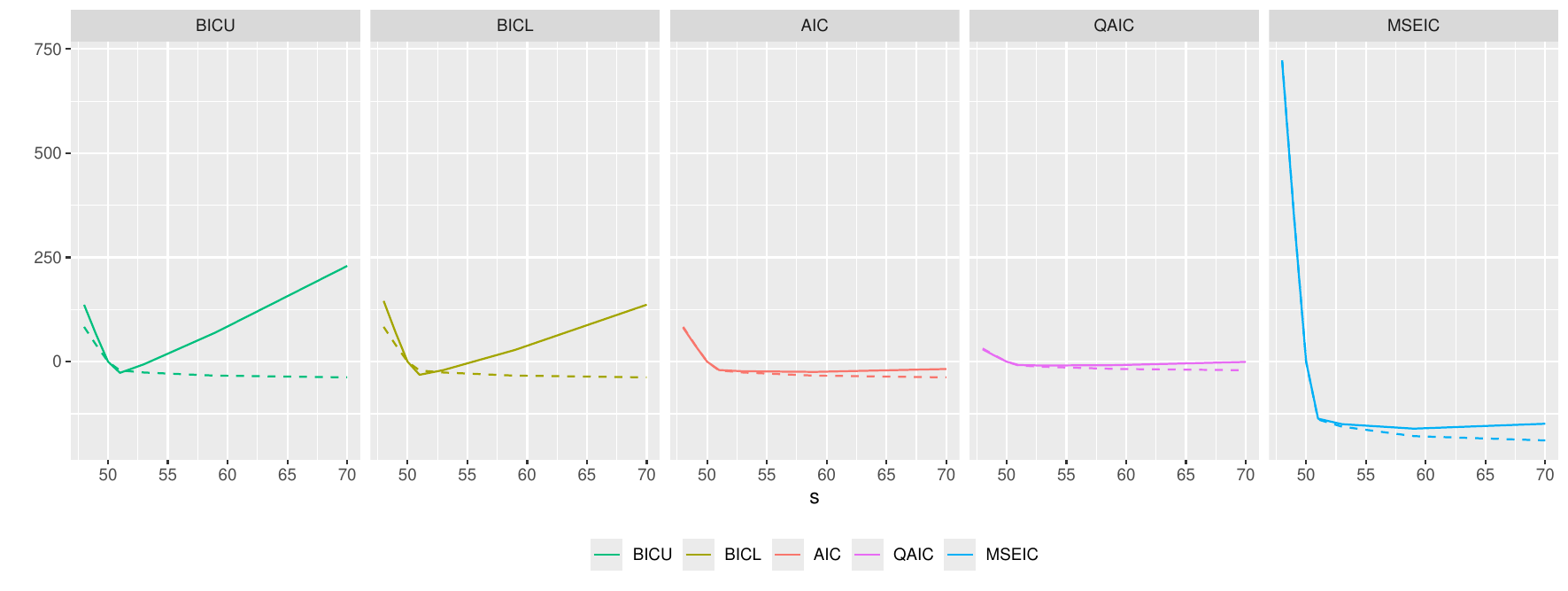}
        
\caption{\footnotesize \textit{
A realization of $\BICU$, $\BICL$, $\AIC$, $\QAIC$ and $\MSEIC$ in the asymptotically dependent model with $ s^*  = 50$, $n = 50.000$ and $ d_n=300$: The solid lines indicate $\mbox{IC}_{k_n}(s)-\mbox{IC}_{k_n}(s^*)$ for $48\leq s\leq 70$ and the dashed lines indicate $[(\mbox{IC}_{k_n}(s)-P(s))-(\mbox{IC}_{k_n}(s^*)-P(s^*)]$ (without the penalty term).
}}
\label{fig:asymp_dep_penalty_figures}
\end{figure}

Next, we present the results of an additional simulation study for an asymptotically dependent model that can also be found in \citet{meyer_tail}, where also directions with $\vb > 1$ are relevant. Let $\bXn$ be an $\Rdn$ valued random vector and $s^*_1, \,s^*_2,\,s^*_3 \in \N \cup \{0 \}$, such that
\begin{align*}
d_n \geq s^*_1 + 2 s^*_2 + 3 s^*_3.
\end{align*}
The parameters $s^*_1,\, s^*_2,\, s^*_3$ specify the number of one-, two-, and three-dimensional extremal directions. In the following, we denote by $\text{Exp}(1)$ the exponential distribution with parameter $1$. The marginal distributions of $\bXn$ are defined by
\begin{align*}
X^{(n)}_j & \sim \text{Pareto}(1), \qquad j = 1, \ldots, s^*_1, \\
(X^{(n)}_j, X^{(n)}_{j+1}) &\sim ( \text{Pareto}(1), X^{(n)}_j + \text{Exp}(1)), \, j = s^*_1 + 1, s^*_1 + 3, \ldots, s^*_1 + 2 \cdot s^*_2 -1, \\
(X^{(n)}_j,  X^{(n)}_{j+1},  X^{(n)}_{j+2}) & \sim ( \text{Pareto}(1), X^{(n)}_j + \text{Exp}(1), X^{(n)}_j + \text{Exp}(1)), \\ & \qquad j = s^*_1 + 2 \cdot s^*_2 + 1, s^*_1 + 2 \cdot s^*_2 + 4, \ldots, s^*_1 + 2 \cdot s^*_2 + 3  \cdot s^*_3 - 2,\\
X^{(n)}_j & \sim \text{Exp}(1), \qquad j =  s^*_1 + 2 \cdot s^*_2 + 3 \cdot  s^*_3 +1, \ldots, d_n.
\end{align*}
In total, there are $\sstar = s^*_1 + s^*_2 + s^*_3$ extremal directions, and the objective is again to identify these directions once more. For the simulation study, we set $s^*_1 = 30, s^*_2 = 15, s^*_3 = 5$ resulting in $s^* = 50$ extremal directions  where $d_n=300$.  Furthermore, we used $n=10.000,25.000$ and $50.000$ with $k_{10.000} = 750$, $k_{25.000}=2.500$  and $k_{50.000}=5.000$, respectively, corresponding to $7.5\%$ - $10\%$ of the data. As before, we chose $k_n$ in such a way that the ratio $\hats/k_n$ remains almost constant, as is confirmed in Figure \ref{fig:asymp_dep_figures}a. In Figure \ref{fig:asymp_dep_figures}b, the boxplots for the empirical estimator for $q$ and in Figure \ref{fig:asymp_dep_figures}c
the boxplots for the empirical estimators for $ g_{\AIC}(q,\mu),g_{\MSEIC}(q,\mu)$ and $g_{\QAIC}(q)$, respectively are presented. These plots resemble the asymptotically independent case in \Cref{fig:asymp_indep_figures}.
Finally, in  \Cref{fig:asymp_dep_directions}, the boxplots of the estimated number of  extremal directions against the sample size on the $x$-axis, and similarly in \Cref{fig:asymp_dep_hellinger}, the boxplots for the Hellinger distance are plotted. In both figures, from left to right, the number of observations increases from $n = 10.000, 25.000$ to $n = 50.000$.  Since \Cref{fig:asymp_dep} and \Cref{fig:asymp_dep_penalty_figures} as well follow a similar style to the asymptotic independent case in \Cref{fig:asymp_indep} and \Cref{fig:asymp_indep_penalty_figures}, respectively, we omit a detailed discussion for brevity. 

\section{Conclusion} \label{sec:conclusion}
In this paper, we extended the four information criteria AIC, BIC, QAIC, and MSEIC to estimate the number of extremal directions from the fixed-dimensional case in \citet{butsch2024fasen} to the high-dimensional case, where the number of bias directions tends to infinity at a similar rate to the number of extremes $k_n$. We further analyzed them for consistency. In contrast to the fixed-dimensional case, where only the Bayesian information criteria and the $\QAIC$ are consistent, we were able to 
provide sufficient criteria for consistency for all information criteria as well as necessary conditions for the $\AIC$ and $\QAIC$. The Bayesian information criteria are consistent without any additional assumption, in analogy with the fixed-dimensional case, except for our fundamental model \Cref{asu:T}.

The simulation studies suggest that the BIC perform very well for large sample sizes; for low sample sizes, however,  underestimation is evident. In contrast, the $\AIC$, $\QAIC$ and $\MSEIC$ perform much better for small sample sizes, but overestimate for large sample sizes. 
In the majority of simulations, we found that the estimates of $\AIC$ and $\MSEIC$ coincide with the number of directions with more than one observation, which explains the overestimation. However, investigating the sample path behaviour of both information criteria reveals an evident change point in either $s^*$ or $s^*+1$. Therefore, these plots of the information criteria are an important graphical tool for the model fit. The reason for detecting a change point at $s^*+1$ rather than at $s^*$ is that, due to the properties of the Euclidean projection, the direction $(1,\ldots,1)$ is estimated as an extremal direction even though it is not. Therefore, the Bayesian criteria also estimate $76$ directions occasionally.

Finally, we would like to point out that the choice of $k_n$ influences the information criteria. Further research into the optimal choice of $k_n$, as in \citet{butsch2024fasen} and \citet{meyer_muscle23}, is needed.

\appendix

\section{Proofs }  \label{proof:information_criteria}

\subsection{Proof of Theorem \ref{th:BIC_Consistency}}~  \label{proof:BIC} $\mbox{}$\\
 (a) \,   \underline{\textbf{Case 1:} }Suppose $s < {\sstar}$. Using the definition of the $\BICU$ gives that 
\begin{align}
& \frac{\BICU_{k_n}(s)  - \BICU_{k_n}({\sstar})}{k_n} \nonumber\\
&\quad = 2\sum_{j=s+1}^{{\sstar}}  \frac{\Tn{j}(k_n)}{k_n} \log \left( \frac{\Tn{j}(k_n)}{k_n} \right) -  2\log \left( \frac{1}{\hats-s} \sum_{j=s+1}^{\hats} \frac{\Tn{j}(k_n)}{k_n} \right) \sum_{i=s+1}^{{\hats}} \frac{\Tn{i}(k_n)}{k_n} \nonumber\\
&\qquad  + 2 \log \left( \frac{1}{{\hats}-{\sstar}} \sum_{j={\sstar}+1}^{\hats} \frac{\Tn{j}(k_n)}{k_n} \right) \sum_{i={\sstar}+1}^{{\hats}} \frac{\Tn{i}(k_n)}{k_n}   \nonumber\\ 
&\qquad + \left\{ 2 \frac{(s - {\sstar}) \log(k_n)}{k_n} + \frac{s}{k_n} \log \left(\frac{ \hats}{2 \pi ( \hats-s) }  \right)   - \frac{\sstar}{k_n} \log \left(\frac{ \hats}{2 \pi ( \hats-{s^*}) }  \right) \right\}\nonumber\\
&\quad = 2 \sum_{j=s+1}^{{\sstar}}  \frac{\Tn{j}(k_n)}{k_n} \log \left( \frac{\Tn{j}(k_n)}{k_n} \right) -  2 \log \left( \frac{1}{{\hats}-s} \sum_{j=s+1}^{\hats} \frac{\Tn{j}(k_n)}{k_n} \right) \sum_{i=s+1}^{\sstar} \frac{\Tn{i}(k_n)}{k_n} \nonumber\\
&\qquad  + 2 \log \left( \frac{{\hats}-s}{{\hats}-{\sstar}} \frac{\sum_{j={\sstar}+1}^{\hats} \Tn{j}(k_n) }{ \sum_{i=s+1}^{\hats} \Tn{i}(k_n)}\right) \sum_{i={\sstar}+1}^{{\hats}} \frac{\Tn{i}(k_n)}{k_n}  \nonumber\\ 
&\qquad  + \left\{2 \frac{(s - {\sstar}) \log(k_n)}{k_n} + \frac{s}{k_n} \log \left(\frac{ \hats}{2 \pi ( \hats-s) }  \right)   - \frac{\sstar}{k_n} \log \left(\frac{ \hats}{2 \pi ( \hats-{s^*}) }  \right)  \right\}\nonumber\\
&\quad =:2\cdot E_{n,1} + 2\cdot E_{n,2} + 2\cdot E_{n,3} + E_{n,4} . \label{eq:AIC2_Cons_pos2} 
\end{align}
Next, we derive the asymptotic behavior of all terms in \Cref{eq:AIC2_Cons_pos2}.\\
\textbf{Step 1.} $E_{n,1}$: A conclusion of Assumption \ref{asu:aic_hd_2} is  as $\ninf$ that
\begin{align}
    E_{n,1} = \sum_{j=s+1}^{{\sstar}}  \frac{\Tn{j}(k_n)}{k_n} \log \left( \frac{\Tn{j}(k_n)}{k_n} \right) \Pconv \sum_{j=s+1}^{{\sstar}} p_j \log \left(  p_j \right). \label{eq:BIC_2_A}
\end{align}
\textbf{Step 2.} $E_{n,2}$: By Assumptions \ref{asu:aic_hd_1}, \ref{asu:aic_hd_2} and  \ref{asu:aic_hd_3}  we have  as $\ninf$ that 
\begin{align}
   E_{n,2} =  -  \log \left( \frac{1}{{\hats}-s} \sum_{j=s+1}^{\hats} \frac{\Tn{j}(k_n)}{k_n} \right) \sum_{i=s+1}^{\sstar} \frac{\Tn{i}(k_n)}{k_n} \Pconv \infty. \label{eq:BIC_2_B}
\end{align}
\textbf{Step 3.} $E_{n,3}$: Similarly, in the case $ c > 0$ a consequence of Assumptions \ref{asu:aic_hd_1}, \ref{asu:aic_hd_2} and  \ref{asu:aic_hd_3}  as $\ninf$ is that
\begin{align}
    E_{n,3} &\centeredop{=} \log \left( \frac{{\hats}-s}{{\hats}-{\sstar}} \frac{\sum_{j={\sstar}+1}^{\hats} \Tn{j}(k_n) }{ \sum_{i=s+1}^{\hats} \Tn{i}(k_n)}\right) \sum_{i={\sstar}+1}^{{\hats}} \frac{\Tn{i}(k_n)}{k_n} \nonumber \\
    &\Pconv \log \left( \frac{\mu }{ c^{-1} \sum_{i=s+1}^{\sstar} p_i + \mu}\right) c \mu. \label{eq:BIC_2_C}
\end{align} 
In the case $c = 0$ we receive by similar arguments and $x \log(x) \to 0$ for $ x \searrow 0$  as $\ninf$ that
\begin{align}
    E_{n,3}     &\Pconv 0.   \label{eq:BIC_2_D} 
\end{align}
\textbf{Step 4.} $E_{n,4}$: Finally, for the last term $E_{n,4}$ we have as $\ninf$,
\begin{align}
    E_{n,4} = 2 \frac{(s - {\sstar}) \log(k_n)}{k_n} + \frac{s}{k_n} \log \left(\frac{ \hats}{2 \pi ( \hats-s) }  \right)   - \frac{\sstar}{k_n} \log \left(\frac{ \hats}{2 \pi ( \hats-{s^*}) }  \right) \Pconv 0.   \label{eq:BIC_2_E}
\end{align}
In summary, \Cref{eq:AIC2_Cons_pos2,eq:BIC_2_A,eq:BIC_2_B,eq:BIC_2_C,eq:BIC_2_D,eq:BIC_2_E} result  as $\ninf$ in
\begin{align*}
    \frac{\BICU_{k_n}(s)  - \BICU_{k_n}({\sstar})}{k_n}  \Pconv \infty.
\end{align*}
\underline{\textbf{Case 2:}} Suppose ${\sstar} < \sn < \pn$ and $\pn/ \sqrt{\hats} = o_\P(1)$.  Again, by using the definition of the $\BICU$, it follows that 
\begin{align}
\BICU_{k_n}&(\sn)  -  \BICU_{k_n}(\sstar) \nonumber\\
&= 2 \sum_{j=\sstar +1}^\sn \Tn{j}(k_n) \left\{ - \log \Big( \frac{\Tn{j}(k_n)}{k_n} \Big) + \log \left( \frac{1}{{\hats}-{\sstar}} \sum_{i={\sstar}+1}^{\hats} \frac{\Tn{i}(k_n)}{k_n} \right)    \right\}\nonumber\\
&\quad + 2\log \left( \frac{\hats - \sn }{\hats - \sstar}  \frac{\sum_{j=\sstar+1}^{\hats} \Tn{j}(k_n)}{\sum_{i=\sn+1}^{\hats} \Tn{i}(k_n)} \right) \sum_{i={\sn}+1}^{\hats} \Tn{i}(k_n)  + 2(\sn - {\sstar}) \log(k_n)  \nonumber\\ 
& \quad + \sn \log \left(\frac{ \hats}{2 \pi ( \hats-\sn) }  \right)   - {s^*} \log \left(\frac{ \hats}{2 \pi ( \hats-{s^*}) }  \right). \label{eq:BIC1_Difference2}
\end{align}
A Taylor expansion of the logarithm around $1$ combined with 
     \begin{align*} 
       \frac{\sum_{j=\sstar+1}^{\hats} \Tn{j}(k_n)}{\sum_{j=\sn+1}^{\hats} \Tn{j}(k_n)} \Pconv 1,
    \end{align*}
 which follows from Assumption \ref{asu:aic_hd_3} and $\sn/\hats \Pconv 0$ as $\ninf$, result in 
 \begin{align}
\BICU_{k_n}&(\sn)  -  \BICU_{k_n}(\sstar) \nonumber\\
&=  2 \sum_{j=\sstar +1}^\sn \Tn{j}(k_n) \left\{ - \log \Big( \frac{\Tn{j}(k_n)}{k_n}   \Big) + \log \left( \frac{1}{{\hats}-{\sstar}} \sum_{i={\sstar}+1}^{\hats} \frac{\Tn{i}(k_n)}{k_n} \right)   \right\} \nonumber\\
&\quad +  2 \frac{\hats - \sn }{\hats - \sstar}  \sum_{j=\sstar+1}^{\hats} \Tn{j}(k_n)  - 2 \sum_{i={\sn}+1}^{\hats} \Tn{i}(k_n)  + 2(\sn - {\sstar}) \log(k_n) \nonumber\\
& \quad + \sn \log \left(\frac{ \hats}{2 \pi ( \hats-\sn) }  \right)   - {s^*} \log \left(\frac{ \hats}{2 \pi ( \hats-{s^*}) }  \right) + o_\P(1) \nonumber\\
&= 2 \sum_{j=\sstar +1}^\sn \Tn{j}(k_n)  \log  \left( \frac{1}{{\hats}-{\sstar}} \sum_{i={\sstar}+1}^{\hats} \frac{\Tn{i}(k_n)}{\Tn{j}(k_n)} \right)  \nonumber\\
&\quad + 2 \frac{(\hats - \sstar) + (\sstar - \sn) }{\hats - \sstar}  \sum_{j=\sstar+1}^{\sn} \Tn{j}(k_n)  +  2\frac{\sstar - \sn }{\hats - \sstar}  \sum_{j=\sn+1}^{\hats} \Tn{j}(k_n) \nonumber\\
& \quad + 2 (\sn - {\sstar}) \log(k_n) + \sn \log \left(\frac{ \hats}{2 \pi ( \hats-\sn) }  \right)   - {s^*} \log \left(\frac{ \hats}{2 \pi ( \hats-{s^*}) }  \right) + o_\P(1) \nonumber\\
&= 2 \sum_{j=\sstar +1}^\sn \Tn{j}(k_n) \left\{ \log  \left( \frac{1}{{\hats}-{\sstar}} \sum_{i={\sstar}+1}^{\hats} \frac{\Tn{i}(k_n)}{\Tn{j}(k_n)} \right) + 1 \right\} \nonumber \\
&\quad +   2 \frac{\sstar- \sn }{\hats - \sstar}   \sum_{j=\sstar+1}^{\sn} \Tn{j}(k_n)  + 2 \frac{\sstar - \sn }{\hats - \sstar}  \sum_{j=\sn+1}^{\hats} \Tn{j}(k_n) + 2(\sn - {\sstar}) \log(k_n) \nonumber\\
& \quad + \sn \log \left(\frac{ \hats}{2 \pi ( \hats-\sn) }  \right)   - {s^*} \log \left(\frac{ \hats}{2 \pi ( \hats-{s^*}) }  \right) + o_\P(1). \label{eq:bic_A}
\end{align} 
Note that the error term of the Taylor expansion of the logarithm is indeed $o_\P(1)$ since $\sn^2 (\Tn{\sstar}(k_n))^2$ $/\hats$ $\Pconv 0$ for $\sstar +1 \le \sn \le q_n$.
By Assumptions \ref{asu:aic_hd_3} and \ref{asu:aic_hd_4} it follows for the first term in \Cref{eq:bic_A} and $\sstar +1 \le j_n \le q_n$ that
\begin{align}
    \bigg| \Tn{j_n}(k_n)  & \log  \left( \frac{1}{{\hats}-{\sstar}} \sum_{i={\sstar}+1}^{\hats} \frac{\Tn{i}(k_n)}{\Tn{j_n}(k_n)} \right) -  \Tn{j_n}(k_n)  \log  \left(  \frac{\mu}{\Tn{j_n}(k_n)} \right) \bigg| \nonumber \\
    &\le \bigg\vert \Tn{\sstar +1}(k_n)   \log \bigg( \frac{1}{{\hats}-{\sstar}} \sum_{i={\sstar}+1}^{\hats} \frac{\Tn{i}(k_n)}{\mu} \bigg) \bigg\vert = o_\P(1). \label{eq:bic_B}
\end{align} 
Moreover, Assumption \ref{asu:aic_hd_4} and $\sstar < \sn \le \pn, \pn/\sqrt{\hats} = o_\P(1)$ yield for the second term in \Cref{eq:bic_A} that
\begin{align}
   \left\vert  \frac{\sstar- \sn }{\hats - \sstar}   \sum_{j=\sstar+1}^{\sn} \Tn{j}(k_n) \right\vert \le    \frac{({\sn} - \sstar)^2}{\hats - \sstar}    \Tn{\sstar+1}(k_n) \le    \frac{({\pn} - \sstar)^2}{\hats - \sstar}    \Tn{\sstar+1}(k_n)  \Pconv 0.  \label{eq:bic_C}
\end{align}
Inserting \Cref{eq:bic_B} and \Cref{eq:bic_C} in \Cref{eq:bic_A} gives
 \begin{align*}
&\BICU_{k_n}(\sn)  -  \BICU_{k_n}(\sstar) \\
&\quad= 2 \sum_{j=\sstar +1}^\sn \left\{ \Tn{j}(k_n) \left( \log  \left(  \frac{\mu}{\Tn{j}(k_n)} \right) + 1 \right)  - \frac{1}{\hats - \sstar}  \sum_{i=\sn+1}^{\hats} \Tn{i}(k_n) \right\} \nonumber\\
& \qquad + 2 (\sn - {\sstar}) \log(k_n) + \sn \log \left(\frac{ \hats}{2 \pi ( \hats-\sn) }  \right)   - {s^*} \log \left(\frac{ \hats}{2 \pi ( \hats-{s^*}) }  \right) + o_\P(1).
\end{align*} 
 Since Assumptions \ref{asu:aic_hd_3} and \ref{asu:aic_hd_4} imply that 
 \begin{align*}
     \sum_{j=\sstar +1}^\sn \left\{ \Tn{j}(k_n) \left( \log  \left(  \frac{\mu}{\Tn{j}(k_n)} \right) + 1 \right)  - \frac{1}{\hats - \sstar}  \sum_{i=\sn+1}^{\hats} \Tn{i}(k_n) \right\} = O_\P(\sn), 
 \end{align*}
we obtain
 \begin{align*}
  & \frac{\BICU_{k_n}(\sn)-  \BICU_{k_n}(\sstar)}{(\sn-\sstar) \log(k_n)} \nonumber \\
&\quad \centeredop{=} \frac{2}{(\sn-\sstar) \log(k_n)} \nonumber \\
&\quad \qquad \cdot \sum_{j=\sstar +1}^\sn  \Bigg\{ \Tn{j}(k_n) \left( \log  \left(  \frac{\mu}{\Tn{j}(k_n)} \right) + 1 \right)  - \frac{1}{\hats - \sstar}  \sum_{i=\sn+1}^{\hats} \Tn{i}(k_n) \Bigg\} \nonumber \\
 & \quad \qquad+ 2  +  \frac{\sn}{(\sn-\sstar) \log(k_n)} \log \left(\frac{ \hats}{2 \pi ( \hats-\sn) }  \right)   - \frac{\sstar}{(\sn-\sstar) \log(k_n)} \log \left(\frac{ \hats}{2 \pi ( \hats-{s^*}) }  \right)\nonumber \\
 &\quad\qquad + o_\P(\log(k_n)^{-1}) \nonumber \\
&\quad\; \Pconv 2  > 0.  
\end{align*} 
(b) \, Note that
\begin{align*}
\BICL_{k_n}(s) = \BICU_{k_n}(s)  - s \log \left( k_n  \right) + s \log \left(\frac{ k_n}{2 \pi \Tn{1}(k_n)}  \right)   - s \log \left(\frac{ \hats}{2 \pi ( \hats-s) } \right).
\end{align*}
By a calculation analog to part (a), the $\BICL$ is also weakly consistent since 
\begin{align*}
\hspace{33pt} \log \left(\frac{ k_n}{2 \pi \Tn{1}(k_n)}  \right)   - \log \left(\frac{ \hats}{2 \pi ( \hats-s) } \right) \Pconv \log \left(\frac{ 1}{ p_{1}}  \right) > 0 \quad \text{ as } n\to\infty. \hspace{33pt} \Box
\end{align*}


\subsection{Proof of Theorem \ref{th:AIC_Consistency}}~  \label{proof:AIC}
\\ \underline{\textbf{Case 1:}} Suppose $s < s^*$. Here we have due to $\log(k_n)/k_n\to 0$ and 
\begin{align*}
\AIC_{k_n}(s) = \frac{1}{2} \BICU_{k_n}(s) + s - s \log \left( k_n  \right) - \frac{s}{2}  \log \left(\frac{ \hats}{2 \pi ( \hats -s) }  \right).
\end{align*}
that
\begin{align*}
    &\frac{\AIC_{k_n}(s)  -  \AIC_{k_n}(s^*)}{k_n} \\
    &\quad=   \frac{\BICU_{k_n}(s) -  \BICU_{k_n}(s^*)}{2 k_n} +  \frac{   s -   s \log \left( k_n  \right) - s \log \left(\frac{ \hats}{2 \pi ( \hats-s) }  \right)}{k_n} \nonumber \\
&\qquad  - \frac{ {s^*} - {s^*} \log \left( k_n  \right) - {s^*} \log \left(\frac{ \hats}{2 \pi ( \hats-{s^*}) }  \right)}{k_n}  \nonumber\\
&\quad \Pconv \infty, 
\end{align*}
where we applied the results of the proof of \Cref{th:BIC_Consistency}, 
and thus, the assertion follows.\\
\underline{\textbf{Case 2:}} Suppose $s^* < \sn < q_n$ and $\pn/ \sqrt{\hats} = o_\P(1)$.  We receive similarly to the proof of \Cref{th:BIC_Consistency} that
 \begin{align*}
 & \AIC_{k_n}(\sn)  -  \AIC_{k_n}(\sstar)  \\
  &\quad =-  \sum_{j=\sstar+1}^{\sn} \Tn{j}(k_n) \log \Big( \frac{\Tn{j}(k_n)}{k_n} \Big) -  \log \Big( \frac{1}{{\hats}-\sn} \sum_{i=\sn+1}^{\hats} \frac{\Tn{i}(k_n)}{k_n} \Big) \sum_{i=\sn+1}^{\hats} \Tn{i}(k_n) \nonumber\\
 &\quad \qquad + \log \bigg( \frac{1}{{\hats}-{\sstar}} \sum_{j={\sstar}+1}^{\hats} \frac{\Tn{j}(k_n)}{k_n} \bigg) \sum_{i={\sstar}+1}^{\hats} \Tn{i}(k_n)   + (\sn - {\sstar}) \\  %
 &\quad  =  \sum_{j=\sstar +1}^\sn \left\{ \Tn{j}(k_n) \left( \log  \left(  \frac{\mu}{\Tn{j}(k_n)} \right) + 1 \right)  - \frac{1}{\hats - \sstar}  \sum_{i=\sn+1}^{\hats} \Tn{i}(k_n) \right\} \nonumber\\
& \quad \qquad +  (\sn - {\sstar})  + o_\P(1).
\end{align*} 
Since the function $f(x) = x \left( \log  \left(  \mu/{x} \right)+ 1 \right) $ with derivative $f'(x) = \log(\mu/x)$ is decreasing for $x \in (\mu, \infty)$ and $\liminf_{\ninf} \Tn{\pn}(k_n) \ge \mu$ $\P$-a.s. (because otherwise we receive a contradiction to  \ref{asu:aic_hd_3}), an application  of Assumptions \ref{asu:aic_hd_3} and \ref{asu:aic_hd_4} with $\Tn{\sstar+1}(k_n) \ge \Tn{\sstar+2}(k_n) \ge \cdots \ge \Tn{\sn}(k_n)$ yields that
 \begin{align*}
  &\frac{\AIC_{k_n}(\sn)  -  \AIC_{k_n}(\sstar)}{\sn - \sstar} \nonumber \\
& \quad \centeredop{\ge}   \biggl\{ \Tn{\sstar+1}(k_n) \biggl( \log  \biggl(  \frac{\mu}{\Tn{\sstar+1}(k_n)} \biggr) + 1 \!  \biggr)  - \frac{1}{\hats - \sstar}  \sum_{i=\sn+1}^{\hats} \! \Tn{i}(k_n) + 1 \biggr\}+ o_\P(1) \nonumber \\
& \quad \Pconv   q \mu \left( \log  \left(  \frac{1}{q} \right) + 1 \right)  - \mu +  1 = \mu \Big( q  \big( 1 - \log  (  q )  \big) - 1  + \frac{1}{\mu} \Big). 
\end{align*}
By assumption, the right-hand side is positive, and the assertion follows.  This condition is also a necessary condition, since if it is not satisfied, we do not have consistency against the $(\sstar +1)$-th model. \hfill$\Box$

\subsection{Proof of Theorem \ref{th:QAIC_Consistency}}~  \label{proof:QAIC} $\mbox{}$\\
\underline{\textbf{Case 1:}} Suppose $s < {\sstar}$. 
    Applying \ref{asu:aic_hd_1}, \ref{asu:aic_hd_2} and \ref{asu:aic_hd_3} gives 
    \begin{align}
        \hats {\widehat{\rho}^{s}}_n = \frac{\hats}{\hats - s}   \sum_{j=s+1}^{\sstar} \frac{\Tn{j}(k_n)}{k_n} +  \frac{\hats}{k_n} \sum_{j=\sstar+1}^{\hats} \frac{\Tn{j}(k_n)}{\hats - s} 
        \Pconv   \sum_{j=s+1}^{\sstar} p_j + c \mu, \label{eq:QAIC_A}
    \end{align}
     as well as
        \begin{align*} 
        \hats {\widehat{\rho}^{{\sstar}}}_n \Pconv c \mu.
    \end{align*}
Therefore, if $c > 0$ we have as $\ninf$    
    \begin{eqnarray*}
        & & \hspace{-0.8cm} \frac{\QAIC_{k_n}(s)  - \QAIC_{k_n}(\sstar)}{\hats} \\
        &=&   -\sum_{j=s+1}^{\sstar} \frac{\log (\hats) + \log \Big( \frac{\Tn{j}(k_n)}{k_n}\Big)}{\hats}   + \frac{\hats - s}{\hats}    \log ( \hats {\widehat{\rho}^{s}}_n)  - \frac{\hats - \sstar}{\hats}   \log( \hats  {\widehat{\rho}^{{\sstar}}}_n)   + \frac{s - \sstar}{\hats}   \\
        &\Pconv&\;  \log \left( \frac{ \sum_{j=s+1}^{\sstar} p_j + c \mu}{ c \mu}  \right) > 0.
    \end{eqnarray*}
    If $c = 0$, then the difference converges to $\infty$ in probability, since $-\log( \hats {\widehat{\rho}^{{\sstar}}}_n ) \Pconv \infty$.
    
\underline{\textbf{Case 2:}} Suppose  $\sstar < \sn < q_n$ and $\pn/ \sqrt{\hats} = o_\P(1)$. We have by definition 
     \begin{align} \label{aa1}
       & \QAIC_{k_n}(\sn)   - \QAIC_{k_n}(\sstar) \\
        &\quad =  \sum_{j={\sstar}+1}^{\sn} \log \left( \Tn{j}(k_n) \right) + ( \hats - \sn) \log \left(    \frac{{\widehat{\rho}^{\sn}}_n}{{\widehat{\rho}^{{\sstar}}}_n} \right) - (\sn - \sstar) \log \left(  k_n {\widehat{\rho}^{{\sstar}}}_n\right) +(\sn-{\sstar}). \nonumber
    \end{align} 
 In the following, we derive an alternative representation of the second summand
 by using a Taylor expansion of the logarithm. Therefore, note that
   \begin{align}
       \Bigl(\frac{  {\widehat{\rho}^{\sn}}_n}{ {\widehat{\rho}^{{\sstar}}}_n} - 1 \Bigr) 
     &=   \frac{1}{k_n \widehat{\rho}^{{\sstar}}_n} (k_n  \widehat{\rho}^{{\sn}}_n - k_n \widehat{\rho}^{{\sstar}}_n )  \nonumber \\
      &=  \frac{1}{k_n \widehat{\rho}^{{\sstar}}_n}   \bigg( \Big(  \frac{\sn - \sstar}{(\hats - \sn)(\hats - \sstar)}\Big) \sum_{j=\sn+1}^{\hats}  \Tn{j}(k_n)  - \frac{1}{\hats - \sn} \sum_{j=\sstar+1}^{\sn} \Tn{j}(k_n)   \bigg)  \nonumber \\
       &= O_\P( \sn/\hats)=o_\P(1),  \label{eq:QAIC_5}
    \end{align}    
   which justifies a Taylor expansion, and the behavior of the error term 
\begin{align*}
      ( \hats - \sn)\left(    \frac{{\widehat{\rho}^{\sn}}_n}{{\widehat{\rho}^{{\sstar}}}_n} - 1 \right)^2= O_\P( s_n^2/\hats)=o_\P(1).
\end{align*}
Then a Taylor expansion of the logarithm and   \ref{asu:aic_hd_4} result in 
\begin{align*}
    &( \hats  -  \sn)  \log \left(    \frac{{\widehat{\rho}^{\sn}}_n}{{\widehat{\rho}^{{\sstar}}}_n} \right) \\
    &\quad= ( \hats - \sn)  \left(    \frac{{\widehat{\rho}^{\sn}}_n}{{\widehat{\rho}^{{\sstar}}}_n} - 1 \right)  + o_\P(1)\\
&\quad=  -( \hats - \sn) \frac{\hats - \sstar}{\hats - \sn} \frac{ \sum_{j=\sstar+1}^{\sn} \Tn{j}(k_n)}{\sum_{j=\sstar+1}^{\hats} \Tn{j}(k_n)}  +  ( \hats - \sn) \big( \frac{\hats - \sstar}{\hats - \sn} - 1 \big)  + o_\P(1)\\
&\quad =  - \sum_{j=\sstar+1}^\sn \frac{\Tn{j}(k_n)}{\mu} + \sn - \sstar + o_\P(1).
\end{align*}
Inserting this in \eqref{aa1} and using
  Assumption  \ref{asu:aic_hd_4} gives
    \begin{align*}
       & \QAIC_{k_n}(\sn)   - \QAIC_{k_n}(\sstar) \\
        &\quad =   \sum_{j={\sstar}+1}^\sn \log \big( \Tn{j}(k_n)\big) - \sum_{j=\sstar+1}^\sn \frac{\Tn{j}(k_n)}{\mu} + 2(\sn - \sstar) - (\sn - \sstar) \log \left(  \mu \right) + o_\P(1)\\
        &\quad =   \sum_{j={\sstar}+1}^\sn \Big( \log \left( \frac{\Tn{j}(k_n)}{\mu} \right) -  \frac{\Tn{j}(k_n)}{\mu}  +2 \Big)+ o_\P(1).
    \end{align*}
    Since the function $\log(x) - x$ is monotone decreasing on $\left[1,\infty\right)$, the $\Tn{j}(k_n)$ are also decreasing in $j$ and lower bounded by $1$. Thus, we obtain with \ref{asu:aic_hd_4} that
        \begin{eqnarray*}
        \frac{\QAIC_{k_n}(\sn)  - \QAIC_{k_n}(\sstar)}{\sn}
        &\ge &  \frac{\sn - \sstar }{\sn}\Big( \log \left( \frac{\Tn{\sstar+1}(k_n)}{\mu} \right) -  \frac{\Tn{\sstar+1}(k_n)}{\mu}  +2 \Big) + o_\P(1)\\
        & \ge &  \frac{1 }{\sstar + 1} \Big( \log \left( q \right) -  q  +2 \Big) + o_\P(1),
    \end{eqnarray*}
    which yields the consistency. Note that this condition is also a necessary condition, since if it is not satisfied, we do not have consistency against the $(\sstar +1)$-th model. \hfill$\Box$

\subsection{Proof of Theorem \ref{th:MSE_Consistency}}~  \label{proof:MSEIC}
$\mbox{}$\\
\underline{\textbf{Case 1:}} Suppose $s < {\sstar}$. Then 
   \begin{align} \label{eq:MSEIC_FF}
        &{\frac{\MSEIC_{k_n}(s)-  \MSEIC_{k_n}(\sstar)}{k_n \hats}} \nonumber \\
        &\quad=   2\frac{s - \sstar}{k_n \hats} + \frac{1}{\hats \widehat{\rho}^{s}_n} \sum_{j=s+1}^{\hats} \left(\frac{\Tn{j}(k_n)}{\hats} - \widehat{\rho}^{{s}}_n \right) ^2    - \frac{1}{\hats \widehat{\rho}^{{\sstar}}_n } \sum_{j={\sstar}+1}^{\hats} \left(\frac{\Tn{j}(k_n)}{k_n} - \widehat{\rho}^{{\sstar}}_n \right) ^2\nonumber\\
        &\quad =: E_{n,1} + E_{n,2} - E_{n,3}. 
    \end{align}
Next, we analyze each of the three terms.\\
\textbf{Step 1.} $E_{n,1}$: By $k_n \hats \Pconv \infty$, follows that 
\begin{align}
    E_{n,1} \Pconv 0. \label{eq:MSEIC_D}
\end{align}
\textbf{Step 2.} $E_{n,2}$: An application of \Cref{asu:T}  and \Cref{eq:QAIC_A} give
\begin{align}
 E_{n,2} &\centeredop{=}  \frac{1}{\hats \widehat{\rho}^{{s}}_n }   \sum_{j=s+1}^{\hats} \left(\frac{\Tn{j}(k_n)}{k_n} - \widehat{\rho}^{{s}}_n \right) ^2 \nonumber \\
 &\centeredop{=}  \frac{1}{\hats \widehat{\rho}^{{s}}_n } \left(\sum_{j=s+1}^{\hats} \left( \frac{\Tn{j}(k_n)}{k_n}  \right)^2 - 2 \widehat{\rho}^{{s}}_n \sum_{j=s+1}^{\hats}   \frac{\Tn{j}(k_n)}{k_n} + \hats (\widehat{\rho}^{{s}}_n )^2 \right) \nonumber \\
 &\Pconv \frac{1}{ \sum_{j=s+1}^{\sstar} p_j + c \mu }  \sum_{j=s+1}^{\sstar} p_j^2. \label{eq:MSEIC_A}
\end{align}
\textbf{Step 3.} $E_{n,3}$: Assumption \ref{asu:aic_hd_4} implies that $0 \le \Tn{j}(k_n) \le \Tn{\sstar+1}(k_n) = O_\P(1)$ for $\sstar +1 \le j$, 
such that \Cref{asu:T}  and  \Cref{eq:QAIC_A} result in
\begin{align}
E_{n,3} &\centeredop{=}  \frac{1}{\hats \widehat{\rho}^{{\sstar}}_n }  \sum_{j=\sstar+1}^{\hats} \left(\frac{\Tn{j}(k_n)}{k_n} - \widehat{\rho}^{{\sstar}}_n \right) ^2 \nonumber \\
 &\centeredop{=} \frac{1}{\hats \widehat{\rho}^{{\sstar}}_n }  \sum_{j=\sstar+1}^{\hats} \frac{\Tn{j}(k_n)^2}{k_n^2} - 2 \frac{1}{\hats} \sum_{j=\sstar+1}^{\hats} \frac{\Tn{j}(k_n)}{k_n} + \frac{1}{\hats } (\widehat{\rho}^{{\sstar}}_n)^2 \nonumber\\
 &\Pconv 0. \label{eq:MSEIC_B}
\end{align}
A consequence of \Cref{eq:MSEIC_FF,eq:MSEIC_A,eq:MSEIC_D,eq:MSEIC_B} is then that
   \begin{align*}
        \frac{\MSEIC_{k_n}(s)-  \MSEIC_{k_n}(\sstar)}{k_n \hats} 
        &\Pconv \frac{1}{ \sum_{j=s+1}^{\sstar} p_j + c \mu} \sum_{j=s+1}^{\sstar} p_j^2 > 0.
    \end{align*} 
     \underline{\textbf{Case 2:}} Suppose  $\sstar < \sn < q_n$ and $\pn/ \sqrt{\hats} = o_\P(1)$.  Since $\widehat{\rho}^{\sn}_n \le \widehat{\rho}^{\sstar}_n$, it follows that 
\begin{align}
    &\MSEIC_{k_n}(\sn) -  \MSEIC_{k_n}(\sstar)  \nonumber  \\
   &\quad = 2(\sn - \sstar) + \frac{k_n}{\widehat{\rho}^{\sn}_n} \sum_{j=\sn+1}^{\hats} \left(\frac{\Tn{j}(k_n)}{k_n} - \widehat{\rho}^{{\sn}}_n \right) ^2    - \frac{k_n}{\widehat{\rho}^{{\sstar}}_n} \sum_{j={\sstar}+1}^{\hats} \left(\frac{\Tn{j}(k_n)}{k_n} - \widehat{\rho}^{{\sstar}}_n \right) ^2 \nonumber \\
   &\quad \ge 2(\sn - \sstar) + \frac{1}{k_n \widehat{\rho}^{\sstar}_n} \left( \sum_{j=\sn+1}^{\hats} \left( \Tn{j}(k_n)  - k_n \widehat{\rho}^{{\sn}}_n \right) ^2    -   \sum_{j={\sstar}+1}^{\hats} \left( \Tn{j}(k_n)  - k_n \widehat{\rho}^{{\sstar}}_n \right) ^2 \right). \label{eq:MSEIC_2}
\end{align}
For the second term in \Cref{eq:MSEIC_2}, we have
\begin{align*}
   &\sum_{j=\sn+1}^{\hats}  \left(  \Tn{j}(k_n)  - k_n \widehat{\rho}^{{\sn}}_n \right) ^2    -   \sum_{j={\sstar}+1}^{\hats} \left( \Tn{j}(k_n)  - k_n \widehat{\rho}^{{\sstar}}_n \right) ^2 \\
   &\quad= - \sum_{j={\sstar}+1}^{\sn} \Tn{j}(k_n)^2 + 2 k_n \widehat{\rho}^{{\sn}}_n  \sum_{j={\sstar}+1}^{\sn} \Tn{j}(k_n) + 2 \sum_{j=\sstar+1}^{\hats} \Tn{j}(k_n) k_n ( \widehat{\rho}^{{\sstar}}_n - \widehat{\rho}^{\sn}_n) \\
   &\quad\quad+ (\sstar - \sn) k_n^2 (\widehat{\rho}^{{\sn}}_n)^2 + (\hats - \sstar)  k_n^2 ( (\widehat{\rho}^{\sn}_n)^2 - (\widehat{\rho}^{\sstar}_n)^2)\\
    &\quad= - \sum_{j={\sstar}+1}^{\sn}\left( \Tn{j}(k_n) - k_n \widehat{\rho}^{{\sn}}_n   \right)^2\\
   &\quad\quad + \left\{2 \sum_{j=\sstar+1}^{\hats} \Tn{j}(k_n) k_n ( \widehat{\rho}^{{\sstar}}_n - \widehat{\rho}^{\sn}_n)  + (\hats - \sstar)  k_n^2 ( (\widehat{\rho}^{\sn}_n)^2 - (\widehat{\rho}^{\sstar}_n)^2) \right\}\\
   &\quad\ge - \sum_{j=\sstar+1}^\sn \left( \Tn{\sstar+1}(k_n) -  k_n \widehat{\rho}^{{\sn}}_n \right)^2  \\
   & \quad\quad+ \left\{  2 \sum_{j=\sstar+1}^{\hats} \Tn{j}(k_n) k_n ( \widehat{\rho}^{{\sstar}}_n - \widehat{\rho}^{\sn}_n)+ (\hats - \sstar)  k_n^2 ( (\widehat{\rho}^{\sn}_n)^2 - (\widehat{\rho}^{\sstar}_n)^2)  \right\}\\
   &\quad=  (\sstar-\sn) \left( \Tn{\sstar+1}(k_n) -  k_n \widehat{\rho}^{{\sn}}_n \right)^2  \\
   &\quad\quad+  (\hats - \sstar) k_n  ( \widehat{\rho}^{{\sstar}}_n - \widehat{\rho}^{\sn}_n) \left\{  2 \frac{1}{(\hats - \sstar)} \sum_{j=\sstar+1}^{\hats} \Tn{j}(k_n)   - k_n ( \widehat{\rho}^{{\sstar}}_n + \widehat{\rho}^{\sn}_n)  \right\}\\   
   &\quad =  E_{n,1} + E_{n,2} \cdot E_{n,3},
\end{align*}
where we used that $ 0 <  \Tn{j}(k_n) -  k_n \widehat{\rho}^{{\sn}}_n < \Tn{\sstar +1}(k_n) -  k_n \widehat{\rho}^{{\sn}}_n $ for $j = \sstar +1, \ldots, \sn$. Now we investigate the three different terms.\\
\textbf{Step 1.} $E_{n,1}$: Since $\Tn{\sstar +1}(k_n) -  k_n \widehat{\rho}^{{\sn}}_n \Pconv q \mu - \mu$  by Assumption \ref{asu:aic_hd_3} and Assumption \ref{asu:aic_hd_4}, it follows that
\begin{align}
    \frac{E_{n,1}}{\sn - \sstar} = -  \left( \Tn{\sstar+1}(k_n) -  k_n \widehat{\rho}^{{\sn}}_n \right)^2 \Pconv -  \left( q \mu - \mu \right)^2 =  - \mu^2 \left( q - 1 \right)^2. \label{eq_MSEIC_A}
\end{align}
\textbf{Step 2.} $E_{n,2}$: By Assumption \ref{asu:aic_hd_3} and  $\Tn{j}(k_n) \le \Tn{\sstar+1}(k_n)= O_\P(1)$ for $j > \sstar$ due to  Assumption \ref{asu:aic_hd_4}, we get 
\begin{align}
  \frac{E_{n,2}}{\sn - \sstar} &=  \frac{(\hats - \sstar)k_n}{\sn - \sstar} ( \widehat{\rho}^{{\sstar}}_n - \widehat{\rho}^{\sn}_n) \nonumber \\
    &=  \frac{1}{\sn - \sstar} \sum_{j=\sstar+1}^{\hats} \Tn{j}(k_n) -  \frac{\sn - \sstar }{(\hats - \sn)(\sn - \sstar)} \sum_{j=\sn+1}^{\hats} \Tn{j}(k_n) \nonumber \\
    &= \frac{1}{\sn - \sstar} \sum_{j=\sstar+1}^{\sn} \Tn{j}(k_n) + \frac{\sstar - \sn}{\sn - \sstar } \frac{1}{\hats - \sn} \sum_{j=\sn+1}^{\hats} \Tn{j}(k_n)    \nonumber \\
    &= O_\P(1). \label{eq_MSEIC_1}
\end{align}
\textbf{Step 3.} $E_{n,3}$: Since $k_n \widehat{\rho}^{{\sstar}}_n \Pconv \mu$ and $k_n\widehat{\rho}^{\sn}_n \Pconv \mu$ by Assumption \ref{asu:T} we get 
\begin{align}
     E_{n,3} &= \left( 2 \frac{1}{\hats - \sstar} \sum_{j=\sstar+1}^{\hats} \Tn{j}(k_n)   - k_n ( \widehat{\rho}^{{\sstar}}_n + \widehat{\rho}^{\sn}_n) \right) \nonumber \\
     &= k_n  \widehat{\rho}^{{\sstar}}_n - k_n \widehat{\rho}^{{\sn}}_n  \nonumber\\
     &=  \frac{1 }{\hats - \sstar} \sum_{j=\sstar+1}^{\hats} \Tn{j}(k_n)  -  \frac{1 }{\hats - \sn} \sum_{j=\sn+1}^{\hats} \Tn{j}(k_n)  \nonumber \\    
     &=  \frac{1 }{\hats - \sstar} \sum_{j=\sstar+1}^{\sn} \Tn{j}(k_n) +  \frac{\sstar - \sn }{(\hats - \sstar)(\hats - \sn)}   \sum_{j=\sn+1}^{\hats} \Tn{j}(k_n)  \nonumber \\
     &= O_\P \left( \frac{s_n}{\hats} \right)=o_\P(1). \label{eq_MSEIC_4}
\end{align}
In summary, with $k_n \widehat{\rho}^{{\sstar}}_n \Pconv \mu$ by Assumption \ref{asu:aic_hd_3}, \Cref{eq_MSEIC_A,eq_MSEIC_1,eq_MSEIC_4} we receive that
\begin{align*}
    \hspace{44pt}  \frac{\MSEIC_{k_n}(\sn) -  \MSEIC_{k_n}(\sstar)}{\sn - \sstar} 
        \ge& 2 + \frac{1}{k_n \widehat{\rho}^{\sstar}_n} \frac{E_{n,1} + E_{n,2} \cdot E_{n,3}}{\sn - \sstar} \nonumber\\
    \Pconv&  2   -  \frac{1}{\mu}  \mu^2 \left( q - 1 \right)^2 = 2 -  (q - 1)^2 \mu. \hspace{44pt} \Box
\end{align*}

\putbib
\end{bibunit}
\end{document}